\documentclass[12pt]{article}
\usepackage{amssymb,amsmath,amsfonts,esint}
\usepackage{relsize}

\newtheorem{theorem}{Theorem}[section]

\newtheorem{proposition}[theorem]{Proposition}
\newtheorem{remark}[theorem]{Remark}
\newtheorem{lemma}[theorem]{Lemma}

\newtheorem{claim}[theorem]{Claim}

\newcommand {\Ac}      {{\mathcal A}}

\newcommand {\Cc}      {{\mathcal C}}
\newcommand {\Dc}      {{\mathcal D}}

\newcommand {\Ic}      {{\mathcal I}}

\newcommand {\Kc}      {{\mathcal K}}
\newcommand {\Lc}      {{\mathcal L}}
\newcommand {\Mc}      {{\mathcal M}}
\newcommand {\Nc}      {{\mathcal N}}
\newcommand {\Pc}      {{\mathcal P}}
\newcommand {\Qc}      {{\mathcal Q}}

\newcommand {\Tc}      {{\mathcal T}}

\newcommand {\R}       {{\bf R}}
\newcommand {\N}       {{\bf N}}

\newcommand {\tH}      {\widetilde{H}}

\newcommand {\tQ}      {\widetilde{Q}}

\newcommand {\tP}      {\widetilde{P}}

\newcommand {\RN}      {\R^n}

\newcommand {\SP}      {Sobolev-Poincar\'{e}~}

\newcommand {\LOP}     {L^1_p(\RN)}
\newcommand {\LMOP}    {L_{p}^{m+1}(\RN)}
\newcommand {\LMP}     {L_{p}^{m}(\RN)}
\newcommand {\LMPR}    {L^m_p(\R)}
\newcommand {\LMIR}    {L^m_\infty(\R)}
\newcommand {\LPRN}    {L_p(\RN)}
\newcommand {\CMD}     {C^{m,(d)}(\RN)}
\newcommand {\CMLD}    {C^{m-1,(d)}(\RN)}

\newcommand {\CM}      {C^{m}(\RN)}
\newcommand {\CMON}    {C^{m-1}(\RN)}
\newcommand {\LIPD}    {\Lip(\RN;d)}

\newcommand {\LIPO}    {\Lip(\RN,d_\omega)}

\newcommand {\intl}    {\int\limits}
\newcommand {\emp}     {\emptyset}

\newcommand {\MC}      {\Mc\Cc}
\newcommand {\MQA}     {M_{q,\eta}(\RN)}

\newcommand {\phq}     {\varphi_q}
\newcommand {\dq}      {\delta_q}
\newcommand {\DPQ}     {\Dc_{p,q}(\RN)}
\newcommand {\PMRN}    {\Pc_{m-1}(\RN)}
\newcommand {\PMP}     {\Pc_{m}(\RN)}
\newcommand {\fs}      {f^{\sharp}_{\infty,E}}

\newcommand {\ARN}     {A_1(\RN)}
\newcommand {\VST}     {\vspace*{1mm}}
\newcommand {\smed}    {\mathlarger{\sum}}
\newcommand {\sbig}    {\mathlarger{\mathlarger{\sum}}}

\DeclareMathOperator{\esssup}{ess\,sup}
\DeclareMathOperator{\essinf}{ess\,inf}
\newcommand {\Lip}     {\operatorname{Lip}}
\newcommand {\diam}    {\operatorname{diam}}
\newcommand {\dist}    {\operatorname{dist}}
\newcommand {\supp}    {\operatorname{supp}}
\newcommand {\MP}      {\operatorname{M}}
\newcommand {\bx}      {\hspace{10mm}$\Box$}
\newcommand {\rbx}     {\hspace{10mm}$\vartriangleleft$}

\newcommand {\nn}      {\nonumber}
\newcommand {\rf}[1]    {(\ref{#1})}      
\newcommand {\reff}[1] {\ref{#1}}         
\newcommand{\lbl}[1]      {\label{#1}}       
\newcommand{\be}          {\begin{eqnarray}}
\newcommand{\bel}[1]      {\begin{eqnarray} \label{#1}}
\newcommand{\ee}           {\end{eqnarray}}
\newcommand {\SECT}[2] {\section*{\centerline{\normalsize
{\bf #1}}} \setcounter{section}{#2}
\setcounter{theorem}{0}\setcounter{equation}{0}}
\begin{document}
\parindent 1em
\parskip 0mm
\medskip
\centerline{\large{\bf Lipschitz spaces generated by the
Sobolev-Poincar\'e inequality}}\vspace*{7mm}
\centerline{\large{\bf and extensions of Sobolev functions. I.}}
\vspace*{10mm} \centerline{By~  {\it Pavel Shvartsman}}\vspace*{5 mm}
\centerline {\it Department of Mathematics, Technion - Israel Institute of Technology}\vspace*{2 mm}
\centerline{\it 32000 Haifa, Israel}\vspace*{2 mm}
\centerline{\it e-mail: pshv@math.technion.ac.il}
\vspace*{10 mm}
\renewcommand{\thefootnote}{ }
\footnotetext[1]{{\it\hspace{-6mm}Math Subject
Classification} 46E35}
\begin{abstract} Let $d$ be a metric on $\RN$ and let
$\CMD$ be the space of $C^m$-function on $\RN$ whose
partial derivatives of order $m$ belong to the space $\LIPD$. We show that the homogeneous Sobolev space
$\LMOP,p>n,$ can be represented as a union of $\CMD$-spaces where $d$ belongs to a family of metrics on $\RN$ with certain ``nice'' properties. This enables us in several important cases to give intrinsic characterizations of the restrictions of Sobolev spaces to arbitrary closed subsets of $\RN$. In particular, we generalize the classical Whitney extension theorem for the space $C^m(\RN)$ to the case of the Sobolev space $\LMP$ whenever $m\ge 1$ and $p>n$.
\end{abstract}
\renewcommand{\contentsname}{ }
\tableofcontents
\addtocontents{toc}{{\centerline{\sc{Contents}}}
\vspace*{10mm}\par}
\SECT{1. Introduction.}{1}
\addtocontents{toc}{~~~1. Introduction.\hfill \thepage\par\VST}
\indent
\par Let $p\in[1,\infty]$ and let $\LOP$ be  the homogeneous Sobolev space consisting of all (equivalence classes of) real valued functions $F\in L_{p,loc}(\RN)$ whose distributional partial derivatives of the first order belong to the space $\LPRN$. We equip $\LOP$ with the seminorm
$$
\|F\|_{\LOP}:=\|\nabla F\|_{\LPRN}.
$$
\par In the case where $p=\infty $, the Sobolev space $L_{\infty }^{1}(\RN)$ can be identified with the space $\Lip(\RN)$ of Lipschitz functions on $\RN$. It is  known that the restriction $\Lip(\RN)|_{E}$ of the Lipschitz space $\Lip(\RN)$ to an {\it arbitrary closed set} $E\subset\RN$ coincides with the space $\Lip(E)$ of Lipschitz functions on $E$. See e.g., \cite{McS}.
Furthermore, the classical Whitney extension operator linearly and continuously maps the space $\Lip(E)$ into the space $\Lip(\RN)$ (see e.g., \cite{St}, Chapter 6).
\par In \cite{Sh2} we have shown that an analogous result also holds for all $p$ in the range $n<p<\infty $,   namely, that {\it very same classical linear Whitney extension operator provides an extension to $\LOP$, $p>n$, even an almost optimal extension, of each function on $E$ which is the restriction to $E$ of a function in $\LOP$.} Using this property of the Whitney extension operator, in \cite{Sh2} we give several intrinsic characterizations of the restriction $\LOP|_E$ whenever $p>n$.
\par In the paper under consideration we present alternative and simpler proofs of some of these results and give new characterizations of the trace space $\LOP|_E$. We also explain the above-mentioned phenomenon of the ``universality" of the Whitney extension operator for the scale $\LOP$, $p>n$.
\par Let $d$ be a metric on $\RN$ and let $\LIPD$ be the space of functions on $\RN$ satisfying the Lipschitz condition with respect to the metric $d$. $\LIPD$ is equipped with the standard seminorm
$$
\|F\|_{\LIPD}:=\sup_{x,y\in\RN,\,x\ne y} \,\,\frac{|F(x)-F(y)|}{d(x,y)}\,.
$$
\par We show that for each $p\in(n,\infty)$ the Sobolev spaces $\LOP$ can be represented as a union of {\it Lipschitz spaces} $\LIPD$ where $d$ belongs to a certain family $\Dc$ of metrics on $\RN$:
\bel{S1-UN}
\LOP=\bigcup_{d\in\Dc}\,\,\LIPD.
\ee
\par Before we describe main ideas of our approach
we need to define several notions and fix some notation: Throughout this paper, the word ``cube'' will mean a closed cube in $\RN$ whose sides are parallel to the coordinate axes. We let $Q(c,r)$ denote the cube in $\RN$ centered at $c$ with side length $2r$. Given $\lambda >0$ and a cube $Q$ we let $\lambda Q$ denote the dilation of $Q$ with respect to its center by a factor of $\lambda$. (Thus $\lambda\,Q(c,r)=Q(c,\lambda r)$.) By $\left|A\right|$ we denote the Lebesgue measure of a measurable set $A\subset \RN$.
\par Representation \rf{S1-UN} is based on the following idea: We slightly modify the classical \SP inequality for $\LOP$-functions, $p>n$, in such a way that the modified inequality can be interpreted as a {\it Lipschitz condition} with respect to a certain metric on $\RN$.
\par Let us recall a variant of the \SP inequality for $p>n$. Let $q\in(n,p]$ and let $F\in \LOP$ be a continuous function. Then for every cube $Q\subset \RN$ and every $x,y\in Q$ the following inequality
\bel{S-P}
|F(x)-F(y)|\leq C(n,q) \diam Q\left(\frac{1}{|Q|}\intl_Q \|\nabla F(u)\|^q\,du\right)^{\frac{1}{q}}
\ee
holds. See, e.g. \cite{M}, p. 61, or \cite{MP}, p. 55.
\par In particular, by \rf{S-P}, for every $x,y\in\RN$
\bel{F-XY}
|F(x)-F(y)|\leq  \|x-y\|\left(\frac{1}{|Q_{xy}|}\intl_{Q_{xy}} h^q(u)\,du\right)^{\frac{1}{q}}
\ee
where
$$
Q_{xy}:=Q(x,\|x-y\|)
$$
and $h=C(n,q)\|\nabla F\|$.
\par Let $h\in L_{p,loc}(\RN)$ be an arbitrary non-negative function.  Inequality \rf{F-XY} motivates us to introduce a function
\bel{DQ}
\dq(x,y:h)=\|x-y\| \left(\frac{1}{|Q_{xy}|}\intl_{Q_{xy}} h(u)^q du\right)^{\frac{1}{q}}~,~~~~~x,y\in \RN.
\ee
\par By this inequality, for each $p\in(n,\infty)$ and every $F\in \LOP$ there exists a non-negative function $h\in \LPRN$ such that $\|h\|_{\LPRN}\le C(n,p,q)\|F\|_{\LOP}$ and
\bel{F-SP}
|F(x)-F(y)|\le \dq(x,y:h)~~~\text{for every}~~~x,y\in\RN.
\ee
\par Conversely, let $n<p<\infty$ and let $F$ be a continuous function on $\RN$. Suppose that there exists a non-negative function $h\in\LPRN$ such that \rf{F-SP} holds. Then $F\in \LOP$ and $\|F\|_{\LOP}\le C(n,p,q)\|h\|_{\LPRN}$. See \cite{H} and the sufficiency part of Theorem \reff{MAIN}, Section 3.
\par Thus the Sobolev space $\LOP$ is completely determined by the ``Lipschitz-like'' conditions \rf{F-SP} with respect to the functions $\dq(h)$ generated by functions $h\in \LPRN$. Of course, these observations lead us to the desired representation \rf{S1-UN} provided $\dq(h)$ is a {\it metric} on $\RN$.
\par However, in general, $\dq(h)$ {\it is not a metric.} Moreover, there are examples of non-negative functions $h\in\LPRN$ such that $\dq(h)$ {\it is not equivalent to any metric on $\RN$}. See, e.g., Semmes \cite{Sm2}, p. 586.
\par Nevertheless, we prove that under a certain additional restriction on  $h$, {\it there does exist a metric $d$ on $\RN$ such that $d\sim\dq(h)$}. More specifically, we show that $\dq(h)$ is equivalent to a metric on $\RN$ whenever {\it $h^q$ belongs to the Muckenhoupt class of $A_1$-weights}. See Section 2 for definition and main properties of $A_1$-weights.
\par We express our result in terms of {\it the geodesic distance} $d_q(h)$ associated with the function $\dq(h)$. Given $x,y\in\RN$ this distance is defined by the formula
\bel{GD-M}
d_q(x,y:h):=\inf\,\smed_{i=0}^{m-1}\,\dq(x_i,x_{i+1}:h)
\ee
where the infimum is taken over all finite sequences of points $\{x_0=x,x_1,...,x_m=y\}$ in $\RN$. Notice that $d_q(h)\le\dq(h)$. Furthermore, $\dq(h)$ is equivalent to a metric on $\RN$ if and only if $\dq(h)$ is equivalent to $d_q(h)$. In Section 2 we prove the following
\begin{theorem}\lbl{MAIN-MT} Let $n\le q<\infty$ and let $h\in L_{q,loc}(\RN)$ be a non-negative function such that $h^q\in \ARN$. Then for every $x,y\in\RN$ the following inequality
\bel{W-G}
d_q(x,y:h)\le \dq(x,y:h)\le C\, d_q(x,y:h)
\ee
holds. Here $C$ is a constant depending only on $n,q,$ and the ``norm'' $\|h^q\|_{A_1}$.
\end{theorem}
\begin{remark} {\em David and Semmes \cite{DS} noted that \rf{W-G} holds for $q=n$. Semmes \cite{Sm2} proved Theorem \reff{MAIN-MT} for a certain family of metric spaces (containing $\RN$) and parameters $q$ big enough. We refer the reader to these papers and also to the works
\cite{Bj,KM,Sm1,Sm3,Sm4} and references therein for
various properties of weights $h$ satisfying inequalities \rf{W-G}.\rbx}
\end{remark}
\begin{remark}\lbl{N-Q} {\em If $h\in L_{q,loc}(\RN)$ is a weight such that $h^q\in \ARN$, then for every $s\in(0,q]$ and every $x,y\in\RN$ the following inequality
\bel{EW-G}
\delta_s(x,y:h)\le\dq(x,y:h)\le \|h^q\|_{A_1}^{\frac1q}\,\delta_s(x,y:h)
\ee
holds. See, e.g., \cite{Sm2}, p. 579, or \cite{GR}; see also Remark \reff{N-Q-PR}.
\par Combining these inequalities with Theorem \reff{MAIN-MT} we conclude that for every $q\in[n,\infty)$ and every $s\in(0,q]$
$$
\delta_s(x,y:h)\sim d_s(x,y:h)\sim d_q(x,y:h),~~~~x,y\in\RN,
$$
with constants in the equivalences depending only on $n,q,s,$ and $\|h^q\|_{A_1}$.\rbx}
\end{remark}
\par Let us formulate the first main result of the paper.
\begin{theorem}\lbl{MAIN} Let $n\le q<p<\infty$. There exists a positive constant $\eta$ depending only on $n,p,$ and $q$ such that the following is true:
\par A continuous function  $F\in\LOP$ if and only if there exists a non-negative function $h\in\LPRN$ such that $h^q\in \ARN$, $\|h^q\|_{A_1}\le \eta,$ and for every $x,y\in\RN$ the following inequality
\bel{LIP1}
|F(x)-F(y)|\le d_q(x,y:h)
\ee
holds. Furthermore,
$$
\|F\|_{\LOP}\sim \inf \|h\|_{\LPRN}.
$$
The constants of this equivalence depend only on $n,p,$ and $q$.
\end{theorem}
\par This result shows that for each $q\in[n,p)$ the Sobolev space $\LOP$ can be represented in the form
\bel{S-UDQ}
\LOP=\bigcup_{d\in\,\DPQ}\LIPD
\ee
where
\bel{D-LIP}
\DPQ:=\{d_q(h): h\in\LPRN, h\ge 0, h^q\in \ARN, \|h^q\|_{A_1}\le \eta(n,p,q)\}.
\ee
\par We apply Theorem \reff{MAIN} to the study of the extension and restriction properties of Sobolev functions. In Section 4 we characterize the restrictions of functions from the space $\LOP$ to an arbitrary closed subset $E$ of $\RN$. Recall that, when $p>n$, it follows from the Sobolev embedding theorem that every function $F\in \LOP$ coincides almost everywhere with a {\it continuous} function on $\RN$. This fact enables us {\it to identify each element} $F\in \LOP$, {\it $p>n$, with its unique continuous representative}. This identification, in particular, implies that $F$ has a well-defined restriction to any given subset of $\RN$.
\par As usual, we define the trace space $\LOP|_E$ of all restrictions of $\LOP$-functions to $E$ by
$$
\LOP|_{E}:=\{f:E\to\R:\text{there exists}~~F\in
\LOP\cap C(\RN)~~\text{such that}~~~F|_{E}=f\}.
$$
We equip this space with the standard quotient space seminorm
$$
\|f\|_{\LOP|_E}
:=\inf \{\|F\|_{\LOP}:F\in\LOP\cap C(\RN), F|_{E}=f\}.
$$
As is customary, we refer to $\|F\|_{\LOP|_E}$ as the {\it trace norm of the function $F$} (in $\LOP$).
\par The main result of Section 4, Theorem \reff{CR-LOP}, provides an explicit formula for the order of magnitude of the trace norm of a function $f:E\to\R$ in the space $\LOP|_E$ whenever $p\in(n,\infty)$.
\begin{theorem}\lbl{CR-LOP} Let $n<p<\infty$ and let $f:\rightarrow \R$ be a continuous function defined on a closed set $E\subset \RN$. Then
\bel{CR-L1P}
\|f\|_{L^1_p(\RN)|_E}\sim \left\{\,\,\intl_{\RN}\left(\,\sup_{y,z\in E}\frac{|f(y)-f(z)|}{\|x-y\|+\|x-z\|}\right)^p\,dx
\right\}^{\frac{1}{p}}.
\ee
The constants of this equivalence depend only on $n$ and $p$.
\end{theorem}
\begin{remark} \lbl{R-L1P}{\em  As we have mentioned above, in \cite{Sh2} we show that the classical extension operator  constructed by Whitney in \cite{W1} for the space $C^1(\RN)$ provides an almost optimal extension operator for the trace space $\LOP|_E$ whenever $n<p<\infty$. In this paper we present a proof of the trace criterion \rf{CR-L1P} which {\it does not use the Whitney extension method}. More specifically, this proof uses only Theorem \reff{MAIN} and McShane's trace theorem \cite{McS}. See Remark \reff{ALG-MC}.\rbx}
\end{remark}
\begin{remark} \lbl{VM-E}{\em Theorem \reff{MAIN} and certain special properties of metrics from $\DPQ$ also enable us to explain the above-mentioned phenomenon of the universality of the Whitney extension operator for the scale of the Sobolev spaces $\LOP$, $p>n$.
\par Let us introduce a special family of metrics on $\RN$. Let $\omega$  be a {\it concave non-decreasing continuous} function on $[0,\infty)$ such that $\omega(0)=0$. We refer to $\omega$ as a ``modulus of continuity''. By $\MC$ we denote the family of all ``moduli of continuity''. For every $\omega\in\MC$ the function
$$
d_{\omega}(x,y)=\omega(\|x-y\|),~~~~ x,y\in\RN,
$$
is a metric on $\RN$. The metric space $(\RN,d_\omega)$ is known in the literature as {\it the metric transform of $\RN$ by $\omega$} or {\it the $\omega$-metric transform} of $\RN$.
\par We show that for each metric $d\in\DPQ$ there exists a mapping $\RN\ni x\mapsto\omega_x\in\MC$ such that
\bel{VAR-D}
d(x,y)\sim \omega _x(\|x-y\|)\,, ~~~~x,y\in\RN,
\ee
with constants depending only on $n,p,$ and $q$. This equivalence motivates us to refer to the metric space $(\RN,d)$ as a \textit{variable metric transform} of $\RN$.
\par It is well-known that for each modulus of continuity $\omega\in\MC$ the Whitney extension operator provides an almost optimal extension of every function $f\in\LIPO|_E$ to a function from $\LIPO$. (See, e.g., \cite{St}, 2.2.3.) In other words, the Whitney operator is a universal extension operator  for the scale $\LIPO$ whenever $\omega\in\MC$. Using \rf{VAR-D} we prove the same statement for the scale $\Lip(\RN;d)$ of the Lipschitz spaces with respect to variable metric transforms $d\in\DPQ$.
This and representation \rf{S-UDQ} imply the universality of the Whitney extension operator for the scale $\LOP,\, p>n$.\rbx}
\end{remark}
\medskip
\par In Section 5 we generalize representation \rf{S-UDQ} to the case of the homogeneous Sobolev space $\LMP$ of order $m\ge 1$. We recall that this space consists of all (equivalence classes of) real valued functions $F\in L_{p,loc}(\RN)$, $p\in[1,\infty],$ whose distributional partial derivatives of order $m$ belong to $\LPRN$. The space $\LMP$ is equipped with the seminorm
$$
\|F\|_{\LMP}:=\|\nabla^mF\|_{\LPRN}
$$
where
\bel{N-M}
\nabla^m F(x):=\left(\smed_{|\alpha|= m}(D^\alpha F(x))^2\right)^{\frac{1}{2}},~~~~x\in\RN.
\ee
\par When $p>n$, by the Sobolev embedding theorem, every function $F\in\LMP$ coincides almost everywhere with a $C^{m-1}$-function, so that we can identify each element $F\in \LMP$ with its unique $C^{m-1}$-representative. As usual, given a closed set $E\subset\RN$ we let $\LMP|_E$ denote the trace of the space $\LMP$ to $E$. This space is normed by 
$$
\|f\|_{\LMP|_E}
:=\inf \{\|F\|_{\LMP}:F\in\LMP\cap C^{m-1}(\RN), F|_{E}=f\}.
$$
\par Let $d$ be a metric on $\RN$ and let $\CMD$ be
the space of functions $F\in\CM$ whose partial derivatives of order $m$ are Lipschitz continuous on $\RN$ with respect to $d$. This space is equipped with the seminorm
$$
\|F\|_{\CMD}:=\smed_{|\alpha|=m}\,\|D^\alpha F\|_{\Lip(\RN;d)}.
$$
\par Using the same approach as for the case $m=1$ we prove that
\bel{DM-PR}
\LMP=\bigcup_{d\in\DPQ}\CMLD
\ee
provided $n\le q<p$. See Theorem \reff{SM-MAIN} and Remark \reff{CM-L}.
\par This representation enables us to prove two extension theorems. The first of them is Theorem \reff{JET-S} which
provides a generalization of the classical Whitney extension theorem for $C^m$-jets \cite{W1} to the case of Sobolev spaces $\LMP$ whenever $n<p<\infty$.
\par Let $J=\{P_x\in\PMRN: x\in E\}$ be a polynomial field, i.e., a mapping which to every point $x\in E$ assigns a polynomial $P_x$ on $\RN$ of degree at most $m-1$. Given an $m$-times differentiable function $F$ and a point $x\in \RN$, we let $T^m_x[F]$ denote the Taylor polynomial of $F$ at $x$ of order $m$. We recall that a variant of the above-mentioned Whitney theorem \cite{W1} for $L^m_\infty$-jets states that there exists a function $F\in L^m_\infty(\RN)\cap \CMON$ such that $T^{m-1}_x[F]=P_x$ for all $x\in E$ if and only if the following quantity
$$
\Nc_{m,\infty}(J):=\sbig_{|\alpha|\le\, m-1}\,\,\,
\sup_{y,z\in E,\, y\ne z } \frac{|D^{\alpha}P_y(y)-D^{\alpha}P_z(y)|}{\|y-z\|^{m-|\alpha|}}
$$
is finite. See also Glaeser \cite{G}.
\begin{theorem} \lbl{JET-S} Let $p\in(n,\infty)$, $m\in\N$,  and let $J=\{P_x\in\PMRN: x\in E\}$ be a polynomial field defined on a closed set $E\subset\RN$.
\par There exists a function $F\in\LMP\cap \CMON$  such that $T^{m-1}_x[F]=P_x$ for every $x\in E$ if and only if the following inequality
\bel{IDF}
\Nc_{m,p}(J):=\sbig_{|\alpha|\le m-1}\,\,
\left(\,\,\intl\limits_{\RN}
\left(\,\,\,
\sup_{y,z\in E}\,\, \frac{|D^{\alpha}P_y(y)-D^{\alpha}P_z(y)|}{\|x-y\|^{m-|\alpha|}+
\|x-z\|^{m-|\alpha|}}\right)^pdx\right)^{\frac{1}{p}}<\infty
\ee
holds. Furthermore,
$$
\Nc_{m,p}(J)\sim \inf\left\{\|F\|_{L^m_p(\RN)}:F\in \LMP\cap\CMON,\, T^{m-1}_x[F]=P_x~\text{for every}~x\in E\right\}.
$$
The constants of this equivalence depend only on $n,m$ and $p$.
\end{theorem}
\par Clearly, for $m=1$ the above equivalence coincides with \rf{CR-L1P}. Similar to the case $m=1$ (see Remark \reff{R-L1P}), we show that the extension operator constructed by Whitney \cite{W1} for the space $C^m(\RN)$ provides an almost optimal extension for every polynomial field $\{P_x\in\PMRN: x\in E\}$ on $E$ whenever $n<p<\infty$ and $E$ is an arbitrary closed subset in $\RN$.
\bigskip
\par Our second extension result provides a constructive characterization of the traces of Sobolev $L^m_p$-functions in {\it one dimensional case}.
\par Let $E$ be a closed subset of $\R$. We recall that
Whitney \cite{W2} constructed a continuous linear extension operator $T_E$ for the trace space $\LMIR|_E$. The operator norm of $T_E$ is bounded by a constant depending only on $m$. In the forthcoming paper \cite{Sh-II} we show that {\it very same Whitney extension operator $T_E$ also provides an almost optimal extension of each function $f\in\LMPR|_E$ to a function from $\LMPR$ whenever $1<p<\infty.$}
\par This leads us to an analogue of the Whitney trace criterion \cite{W2} for the space $\LMIR$. Recall that, by this criterion, for arbitrary positive integer $m$ and every function $f$ on $E$,
\bel{T-R1}
\|f\|_{\LMIR|_E}\sim
\sup_{S\subset E,\,\,\# S=m+1}
|\Delta^mf[S]|.
\ee
In this formula given a finite set $S=\{x_0,x_1,...,x_m\}\subset \R$ the quantity $\Delta^mf[S]$ denotes {\it the divided difference} of $f$ on $S$ of order $m$:
$$
\Delta^mf[S]=\Delta^mf[x_0,x_1,...,x_{m}]
=\sbig_{i=0}^m\,\,\frac{f(x_i)}
{\omega_S'(x_i)}
$$
where $\omega_S(x)=(x-x_0)\cdot(x-x_1)\cdot...\cdot(x-x_m)$. In \rf{T-R1} the symbol $\#$ denotes the number of points of a set. (See also J. Merrien \cite{Mer} and Glaeser \cite{G}.)
\par In \cite{Sh-II} we present a counterpart of formula \rf{T-R1}. We prove that for every $1<p<\infty$, every closed set $E\subset\R$ and every function $f$ on $E$
\be
\|f\|_{\LMPR|_E}&\sim& \left\{ \intl\limits_{\R}\sup_{S\subset E,\,\,\# S=m+1}\left(\frac{|\Delta^mf[S]|\diam S}{\diam (\{x\}\cup S)}\right)^pdx\right\}^{\frac{1}{p}}\nn\\\nn\\
&\sim& \left\{ \intl\limits_{\R}\sup_{\substack{ x_0<x_1<...<x_m,\\x_i\in E}}\,\,\frac{|\,\Delta^{m-1}f[x_0,...,x_{m-1}]-
\Delta^{m-1}f[x_1,...,x_{m}]|^{\,p}}{|x-x_0|^p+|x-x_m|^{p}}
\,\,dx\right\}^{\frac{1}{p}}.\nn
\ee
The constants in these equivalences depend only on $p$ and $m$. \medskip
\par As we have noted above, the extension operator $T_E$ constructed by Whitney \cite{W2} for the space $C^m(\R)|_E$ provides an almost optimal extension for every trace space $L^m_p(\R)|_E$ and {\it every} closed set $E\subset\R$ whenever $p>1$. Since $T_E$ is {\it linear}, we obtain the following result.

\begin{theorem} \lbl{LINEXT} For every closet set
$E\subset\R$ and every $p>1$ there exists a linear extension operator which maps the trace space $\LMPR|_E$ continuously into $\LMPR$. Its operator norm is bounded by a constant depending only on $p$.
\end{theorem}
\par Note that G. K. Luli \cite{L} proved the existence of a continuous linear extension operator for the space $L^m_p(\R)|_E,$ $p>1,$ for every {\it finite} set $E\subset\R$. We also remark that an analog of Theorem \reff{LINEXT} for the space $L^2_p(\R^2)$, $p>2$, has been proven in recent works of A. Israel \cite{Is} and the author \cite{Sh3}. Quite recently  C. Fefferman, A. Israel and G. K. Luli \cite{FIL} proved the existence of a continuous linear extension operator for the space $\LMP|_E$ whenever $n<p<\infty$ and $E\subset\RN$ is an arbitrary closed set.

\par {\bf Acknowledgements.} I am very thankful to M. Cwikel for useful suggestions and remarks. The results of this paper were presented at ``Whitney Problems Workshop", Palo-Alto, August 2010. I am very grateful to C. Fefferman, N. Zobin and all participants of this conference for stimulating discussions and valuable advice. I am also pleased to thank P. Haj{\l}asz who kindly drew the author's attention to several important papers related to Theorems \reff{MAIN-MT} and \reff{MAIN}.
\SECT{2. Metrics on $\RN$ generated by the Sobolev-Poincar\'e inequality.}{2}
\addtocontents{toc}{2. Metrics on $\RN$ generated by the Sobolev-Poincar\'e inequality. \hfill\thepage\par\VST}
\indent
\par Let us fix additional notation. Throughout the paper $\gamma,\gamma_1,\gamma_2...,$ and $C,C_1,C_2,...$ will be generic positive constants which depend only on parameters determining function spaces ($p,q,n,m,$ etc). These constants can change even in a single string of estimates. The dependence of a constant on certain parameters is expressed, for example, by the notation $C=C(n,p,q)$. We write $A\sim B$ if there is a constant $C\ge 1$ such that $A/C\le B\le CA$.
\par Given $x=(x_1,x_2,...,x_n)\in\RN$ by $\|x\|:=\max \{|x_1|,|x_2|,...,|x_n|\}$ we denote the uniform norm in $\RN$. Let $A,B\subset \RN$. We put
$\diam A:=\sup\{\|a-a'\|:~a,a'\in A\}$ and
$$
\dist(A,B):=\inf\{\|a-b\|:~a\in A, b\in B\}.
$$
For $x\in \RN$ we also set $\dist(x,A):=\dist(\{x\},A)$.
We put  $\dist(A,B)=+\infty$ and $\dist({x},B)=+\infty$ whenever $B=\emp.$ For each pair of points $z_1$ and $z_2$ in $\RN$ we let $(z_1,z_2)$ denote the open line segment joining them. Finally given a cube $Q$ in $\RN$ by $c_Q$ we denote its center, and by $r_Q$ a half of its side length. (Thus $Q=Q(c_Q,r_Q).$)
\par We let $\PMP$ denote the space of polynomials of degree at most $m$ defined on $\RN$. Finally given a function $g\in L_{1,loc}(\RN)$ we let $\Mc[g]$ denote its Hardy-Littlewood maximal function:
$$
\Mc[g](x):=\sup_{Q\ni x}\frac{1}{|Q|}\int_Q|g|dx,~~~~x\in\RN.
$$
Here the supremum is taken over all cubes $Q$ in $\RN$ containing $x$.
\bigskip
\par {\bf 2.1 Proof of Theorem \reff{MAIN-MT}.}
\addtocontents{toc}{~~~~2.1 Proof of Theorem \reff{MAIN-MT}. \hfill \thepage\par}
Let $w$ be a weight on $\RN$, i.e., a non-negative locally integrable function. Let $q>0$ and let $\phq(w):\RN\times\RN\to\R_+$ be a function defined by the formula
\bel{DF-PHI}
\phq(x,y:w):=\|x-y\| \sup_{Q\ni x,y} \left(\frac{1}{|Q|}\intl_{Q} w(u)\,du\right)^{\frac{1}{q}},~~~x,y\in\RN.
\ee
\begin{proposition}\lbl{P3.1} Let $n\le q<\infty$. Then for every  $x,y\in \RN$ and every finite family of points
$x_0=x,x_1,...x_{m-1},x_m=y$ in $\RN$ the following inequality
\bel{3.4'}
\phq(x,y:w)\leq 16\, \smed^{m-1}_{i=0}\,\phq(x_i,x_{i+1}:w)
\ee
holds.
\end{proposition}
\par \textit{Proof.} Let $K$ be an arbitrary cube in $\RN$ and let $x,y\in K$. Prove that
\bel{K-I}
\|x-y\| \left(\frac{1}{|K|}\intl_{K} w(u)\,du\right)^{\frac{1}{q}}\leq 16\,\smed^{m-1}_{i=0}\,\phq(x_i,x_{i+1}:w).
\ee
Let
$$
\tQ=Q(x,2\|x-y\|)=2\,Q_{xy}.
$$
\par Consider two cases.
\par \textit{The first case.} Suppose that {\it there exists $j\in\{0,...,m-1\}$ such that $x_j\in 2\,Q_{xy}$, but}
$$
x_{j+1}\notin 4\,Q_{xy}.
$$
Hence $\|x-x_j\|\le 2\|x-y\|$ and $\|x-x_{j+1}\|\ge 4\|x-y\|$ so that
\bel{3.6}
\|x_i-x_{j+1}\|\ge 2\|x-y\|.
\ee
Since $x,y\in K$, we have  $\diam K\geq \|x-y\|$ which easily implies the inclusion
\bel{QXY-K}
2\,Q_{xy}\subset 5K.
\ee
\par Hence $x_j\in 2\,Q_{xy}\subset 5K.$\medskip
\par Consider the following two subcases. First let us assume that
$$
x_{j+1}\in 8K.
$$
Since $x_j\in 2\,Q_{xy}\subset 5K$, we have $x_j,x_{j+1}\in 8K$. Hence, by \rf{3.6} and definition \rf{DF-PHI} of $\phq$, we obtain
\be
\|x-y\|\left(\frac{1}{|K|}\intl_Kw(u)du\right) ^{\frac{1}{q}}&\leq& \|x_j-x_{j+1}\|\left(\frac{1}{|K|}\intl_Kw(u)du\right) ^{\frac{1}{q}}\nn\\
&\leq& 8^{\frac{n}{q}}\|x_j-x_{j+1}\|\left(\frac{1}{|8K|} \intl_{8K}w(u)du\right)^{\frac{1}{q}}\nn\\&\leq& 8^{\frac{n}{q}}\phq(x_j,x_{j+1}:w).\nn
\ee
Since $n\le q$, we have
$$
\|x-y\|\left(\frac{1}{|K|}\intl_Kw(u)du\right) ^{\frac{1}{q}}\le 8\phq(x_j,x_{j+1}:w)
$$
proving \rf{K-I} in the case under consideration.
\par Now consider the second subcase where
$$
x_{j+1}\notin 8K.
$$
Since $x_j\in 5K$, we conclude that
$\|x_j-x_{j+1}\|\ge 3r_K.$
\par Let $\overline{Q}:=Q(x_j,2\|x_j-x_{j+1}\|).$ Since $x_j\in 5K$, for every $z\in K$ we have
$$
\|x_j-z\|\le \|x_j-c_K\|+\|c_K-z\|\le 5r_K+r_K=6r_K \le 2\|x_j-x_{j+1}\|=r_{\overline{Q}}
$$
proving that $\overline{Q}\supset K$. Clearly,
$r_K\le r_{\overline{Q}}$ and $\overline{Q}\ni x_j,x_{j+1}$.
\par Hence
\be
\|x-y\|\left(\frac{1}{|K|}\intl_Kw(u)\,du\right)
^{\frac{1}{q}}&=&
2^{\frac{n}{q}}\|x-y\|r_K^{-\frac{n}{q}}
\left(\intl_Kw(u)\,du\right)
^{\frac{1}{q}}\nn\\
&\le& 2^{\frac{n}{q}}\|x-y\|r_K^{-\frac{n}{q}}
\left(\intl_{\overline{Q}} w(u)\,du\right)
^{\frac{1}{q}}.\nn
\ee
Since $x,y\in K$, we have $\|x-y\|\le \diam K=2r_K$. Combining this with inequality $r_K\le r_{\overline{Q}}$\,, we obtain
\be
\|x-y\|\left(\frac{1}{|K|}\intl_K w(u)\,du\right)
^{\frac{1}{q}}&\le&
2^{\frac{n}{q}+1}r_K^{1-\frac{n}{q}}
\left(\intl_{\overline{Q}}
w(u)\,du\right)^{\frac{1}{q}}
\le 2^{\frac{n}{q}+1}r_{\overline{Q}}^{1-\frac{n}{q}}
\left(\intl_{\overline{Q}}
w(u)\,du\right)^{\frac{1}{q}}\nn\\
&=&2^{\frac{2n}{q}+1}r_{\overline{Q}}
\left(\frac{1}{|\overline{Q}|}
\intl_{\overline{Q}}w(u)\,du\right)
^{\frac{1}{q}}.\nn
\ee
Since $r_{\overline{Q}}=2\,\|x_j-x_{j+1}\|$ and $\frac{n}{q}\le 1$, we have
$$
\|x-y\|\left(\frac{1}{|K|}
\intl_{K}w(u)\,du\right)^{\frac{1}{q}}\le 2^4\|x_j-x_{j+1}\|\left(\frac{1}{|\overline{Q}|}
\intl_{\overline{Q}}w(u)\,du\right)
^{\frac{1}{q}}.
$$
But $x_j,x_{j+1}\in \overline{Q}$ so that, by definition \rf{3.4'},
$$
\|x-y\|\left(\frac{1}{|K|}\intl_{K}w(u)\,du\right)
^{\frac{1}{q}}\le 2^4\phq(x_j,x_{j+1}:w)
$$
proving \rf{K-I}.
\par \textit{The second case}. Suppose that the assumption of the first case is not satisfied, i.e., {\it for each $i\in\{0,...,m-1\}$ such that $x_i\in 2\, Q_{xy}$ we have  $x_{i+1}\in 4\, Q_{xy}$}.
\par Let us define a number $j\in \{0,1,...,m-1\}$ as follows. If $\{x_0,x_1,...,x_m\}\subset 2\,Q_{xy},$ we put $j=m-1$. If $\{x_0,x_1,...,x_m\}\nsubseteq 2\,Q_{xy},$ then there exists $j\in {0,1,...,m-1}$, such that
$$
\{x_0,x_1,...,x_j\}\subset 2\,Q_{xy}
$$
but
$$
x_{j+1}\notin 2\,Q_{xy}.
$$
Note that, by the assumption, $x_{j+1}\in 4\,Q_{xy}$.
\par Prove that
$$
\|x-y\|\left(\frac{1}{|K|}\intl_{K}w(u)\,du\right)
^{\frac{1}{q}}\le 16\,\smed_{i=0}^j\,\phq(x_i,x_{i+1}:w).
$$
In fact, since $x_{j+1}\notin 2\,Q_{xy}=Q(x,2\|x-y\|)$ we have
$$
\|x_0-x_{j+1}\|=\|x-x_{j+1}\|\ge 2\|x-y\|
$$
so that
\bel{3.11}
\smed_{i=0}^j\,\|x_i-x_{i+1}\|\ge \|x_0-x_{j+1}\|\ge 2\,\|x-y\|.
\ee
\par Recall that, by \rf{QXY-K}, $2\,Q_{xy}\subset 5K$ so that $10K\supset 4\,Q_{xy}.$ Hence
\bel{3.12}
x_0,x_1,...,x_j,x_{j+1}\in 10K.
\ee
By \rf{3.11},
\be
\|x-y\|\left(\frac{1}{|K|}
\intl_{K}w(u)\,du\right)^{\frac{1}{q}}&\le& \frac{1}{2}\left(\smed_{i=0}^j\,\|x_i-x_{i+1}\|\right)
\left(\frac{1}{|K|}
\intl_{K}w(u)\,du\right)^{\frac{1}{q}}\nn\\
&\le&\frac{(10)^{\frac{n}{q}}}{2}
\left(\smed_{i=0}^j\,\|x_i-x_{i+1}\|\right)
\left(\frac{1}{|10K|}
\intl_{10K}w(u)\,du\right)^{\frac{1}{q}}\nn.
\ee
Since $10^{\frac{n}{q}}\le 10$, we obtain
$$
\|x-y\|\left(\frac{1}{|K|}\intl_{K}w(u)\,du\right)
^{\frac{1}{q}}\le
5\left(\smed_{i=0}^j\,\|x_i-x_{i+1}\|\right)
\left(\frac{1}{|10K|}\intl_{10K}w(u)\,du\right)^{\frac{1}{q}}.
$$
But, by \rf{3.12}, $x_i,x_{i+1}\in 10K$ for every $0\le i\le j$, so that
$$
\|x_i-x_{i+1}\|\left(\frac{1}{|10K|}
\intl_{10K}w(u)\,du\right)^{\frac{1}{q}}\le
\phq(x_{i},x_{i+1}:w).
$$
Hence
$$
\|x-y\|\left(\frac{1}{|K|}
\intl_{K}w(u)\,du\right)^{\frac{1}{q}}\le
5\,\smed_{i=0}^j\,\phq(x_{i},x_{i+1}:w)\le
5\,\smed_{i=0}^{m-1}\,\phq(x_{i},x_{i+1}:w)
$$
proving inequality \rf{K-I}.
\par Thus \rf{K-I} is proven for an arbitrary family of points $x_0=x,x_1,...,x_m=y$ in $\RN$. Taking the supremum in this inequality over all cubes $K\ni x,y$, we finally obtain the statement of the proposition.\bx
\par We recall that a non-negative function $w\in L_{1,loc}(\RN)$ is said to be $A_1$-weight if there exists a constant $\lambda >0$ such that for every cube $Q\subset \RN$ the following inequality
\bel{4.3}
\frac{1}{|Q|}\int_Qw(u)\,du\le \lambda \essinf\limits_Q w
\ee
holds. See, e.g. \cite{GR}. We put $\|w\|_{A_1}=\inf \lambda$.
\par Clearly, a weight $w\in A_1$ if and only if there exists $\lambda>0$ such that
$$
\Mc[w]\le \lambda w(x)~~~~~a.e.~on~ \RN.
$$
Furthermore,
$$
\|w\|_{A_1}\sim \esssup_{\RN}\frac{\Mc[w](x)}{w(x)}.
$$
\par Let us also notice the following important property of $A_1$-weights: for every $w\in \ARN$ and every two cubes $K,Q$ in $\RN$ such that $Q\subset K$  the following inequality
\bel{4.C}
\frac{1}{|K|}\intl_K w(u)du\le \|w\|_{A_1}\frac{1}{|Q|}\intl_Q w(u)du
\ee
holds.
\par These properties of $A_1$-weights and Proposition \reff{P3.1} enable us to finish the proof of Theorem \reff{MAIN-MT} as follows.\medskip
\par The first inequality in \rf{W-G} is trivial so let us prove that for every $x,y\in\RN$
\bel{T-L}
\dq(x,y:h)\le C\, d_q(x,y:h)
\ee
provided $q\in[n,\infty)$ and $w:=h^q\in \ARN$. Here $C$ is a constant depending only on $n,q,$ and $\|w\|_{A_1}$. In fact, by \rf{DQ}, for every $x,y\in \RN$
$$
\dq(x,y:h):=\|x-y\| \left(\frac{1}{|Q_{xy}|}\intl_{Q_{xy}} h(u)^q du\right)^{\frac{1}{q}}=\|x-y\| \left(\frac{1}{|Q_{xy}|}\intl_{Q_{xy}} w\,du\right)^{\frac{1}{q}}
$$
so that, by \rf{DF-PHI},
\bel{D-LPH}
\dq(x,y:h)\le
\|x-y\| \sup_{Q\ni x,y} \left(\frac{1}{|Q|}\intl_{Q} w(u)\,du\right)^{\frac{1}{q}}=\phq(x,y:w).
\ee
\par On the other hand, the cube $Q_{xy}\ni x,y$ and $\diam Q_{xy}=2\|x-y\|$ so that for every cube $K\ni x,y$ we have $Q_{xy}\subset 3K$. Hence, by \rf{4.C},
$$
\frac{1}{|K|}\intl_{K} w\,du\le  \frac{3^n}{|3K|}\intl_{3K} w(u)\,du\le 3^n\|w\|_{A_1}\left(\frac{1}{|Q_{xy}|}\intl_{Q_{xy}} w\,du\right).
$$
Combining this inequality with \rf{DF-PHI}, we obtain
$$
\phq(x,y:w):=\|x-y\| \sup_{K\ni x,y} \left(\frac{1}{|K|}\intl_{K} w\,du\right)^{\frac{1}{q}}\le
3^{\frac{n}{q}}\|w\|_{A_1}^{\frac1q}\,\|x-y\|
\left(\frac{1}{|Q_{xy}|}\intl_{Q_{xy}} w\,du\right)^{\frac{1}{q}}
$$
proving that
\bel{PH-D}
\phq(x,y:w)\le C\dq(x,y:h).
\ee
\par Let $\{x_0=x,x_1,...,x_m=y\}$ be an arbitrary family of points in $\RN$. Then, by \rf{D-LPH} and Proposition \reff{P3.1},
$$
\dq(x,y:h)\le\phq(x,y:w)\le 16\, \smed^{m-1}_{i=0}\, \phq(x_i,x_{i+1}:w)
$$
so that, by \rf{PH-D},
$$
\dq(x,y:h)\le C\,\smed^{m-1}_{i=0}\, \dq(x_i,x_{i+1}:h).
$$
Taking the infimum in this inequality over all families of points $\{x_0=x,x_1,...,x_m=y\}$, we obtain the required inequality \rf{T-L}. See \rf{GD-M}.
\par Theorem \reff{MAIN-MT} is completely proved.\bx
\begin{remark}\lbl{N-Q-PR} {\em Prove inequality \rf{EW-G}. Since $0<s\le q$, the first inequality in \rf{EW-G} is elementary. Prove the second inequality. Let $w:=h^q$. Then, by \rf{4.3}, for every $x,y\in\RN$
$$
\frac{1}{|Q_{xy}|}\int_{Q_{xy}} w(u)\,du\le \|w\|_{A_1} \essinf\limits_{Q_{xy}} w
$$
so that
$$
\left(\frac{1}{|Q_{xy}|}\int_{Q_{xy}} w(u)\,du\right)^{\frac1q}\le \|w\|_{A_1}^{\frac1q} \left(\essinf\limits_{Q_{xy}} w^{\frac{s}{q}}\right)^{\frac1s}
\le \|w\|_{A_1}^{\frac1q}
\left(\frac{1}{|Q_{xy}|}\int_{Q_{xy}} w^{\frac{s}{q}}(u)\,du\right)^{\frac1s}.
$$
Hence
$$
\left(\frac{1}{|Q_{xy}|}\int_{Q_{xy}} h^q(u)\,du\right)^{\frac1q}
\le \|h^q\|_{A_1}^{\frac1q}
\left(\frac{1}{|Q_{xy}|}\int_{Q_{xy}} h^s(u)\,du\right)^{\frac1s}
$$
proving \rf{EW-G}.\rbx}
\end{remark}
\bigskip
\par {\bf 2.2 Variable metric transforms.}
\addtocontents{toc}{~~~~2.2 Variable metric transforms. \hfill \thepage\par\VST}
Let $q\in[n,\infty)$ and let $\eta>0$. Define a family of metrics on $\RN$
$$
\MQA=\{d_q(h):h^q\in A_1(\RN), \|h^q\|_{A_1}\le\eta\}.
$$
\par In this subsection we present several simple but important properties of metrics from the family $\MQA$. Note that, by Theorem \reff{MAIN-MT} and Remark \reff{N-Q}, for every $s\in(0,q]$ and every metric $d=d_q(h)\in\MQA$ such that $h^q\in A_1(\RN), \|h^q\|_{A_1}\le\eta$ we have
\bel{F-DS}
d(x,y)\sim\delta_s(x,y:h)=\|x-y\| \left(\frac{1}{|Q_{xy}|}\intl_{Q_{xy}} h(u)^s du\right)^{\frac{1}{s}},~~~~x,y\in\RN.
\ee
with constants depending only on $n,q,$ and $\eta$.
\begin{claim} Let $n\le s\le q$ and let $d=d_q(h)\in\MQA$ where $h\in A_1(\RN), \|h^q\|_{A_1}\le\eta$. Then for every $x,y\in\RN$ the following inequality
$$
d(x,y)\sim \inf_{Q\ni x,y}\diam Q\left(\frac{1}{|Q|}\intl_{Q}
h^sdu\right)^{\frac{1}{s}}
$$
holds with constants depending only on $n,q,$ and  $\eta$. Here the infimum is taken over all cubes $Q\ni x,y$.
\end{claim}
\par {\it Proof.} Since $x,y\in Q_{xy}=Q(x,\|x-y\|)$ and $\diam Q_{xy}=2\|x-y\|$,
\be
I(x,y)&:=&\inf_{Q\ni x,y}\diam Q\left(\frac{1}{|Q|}\intl_{Q}
h^sdu\right)^{\frac{1}{s}}\le \diam Q_{xy}\left(\frac{1}{|Q_{xy}|}\intl_{Q_{xy}}
h^s\,du\right)^{\frac{1}{s}}\nn\\&\le& 2\|x-y\|\left(\frac{1}{|Q_{xy}|}\intl_{Q_{xy}}
h^s\,du\right)^{\frac{1}{s}}\nn
\ee
so that, by \rf{F-DS},
$$
I(x,y)\le C(n,q) d(x,y).
$$
\par On the other hand, for every cube $Q\subset\RN$ such that $Q\ni x,y$ we have $3Q\supset Q_{xy}$. Hence, by \rf{F-DS},
\be
d(x,y)&\sim&\|x-y\| \left(\frac{1}{|Q_{xy}|}\intl_{Q_{xy}} h(u)^s du\right)^{\frac{1}{s}}\sim \|x-y\|^{1-\frac{n}{s}} \left(\,\intl_{Q_{xy}} h(u)^s du\right)^{\frac{1}{s}}\nn\\&\le&
(\diam Q)^{1-\frac{n}{s}}\left(\,\,\intl_{3Q} h(u)^s du\right)^{\frac{1}{s}}\le C \diam Q\left(\frac{1}{|3Q|}\intl_{3Q} h(u)^s du\right)^{\frac{1}{s}}.\nn
\ee
\par Since $h\in A_1(\RN)$ and $\|h^q\|_{A_1}\le\eta$, by \rf{4.C},
\be
\left(\frac{1}{|3Q|}\intl_{3Q} h(u)^s du\right)^{\frac{1}{s}}&\le&
\left(\frac{1}{|3Q|}\intl_{3Q} h(u)^q du\right)^{\frac{1}{q}}\nn\\
&\le& \eta^{\frac{1}{q}}\essinf_{3Q} h\le \eta^{\frac{1}{q}}\essinf_{Q} h\le \eta^{\frac{1}{q}}\left(\frac{1}{|Q|}\intl_{Q} h(u)^s du\right)^{\frac{1}{s}}.\nn
\ee
Hence
$$
d(x,y)\le C\diam Q\left(\frac{1}{|Q|}\intl_{Q}
h^sdu\right)^{\frac{1}{s}}
$$
proving the claim.\bx
\begin{claim}\lbl{VMT} Every metric $d\in\MQA$ has the following properties:\medskip
\par (a). There exists a mapping $\RN\in x\to \omega_x\in \MC$ such that
\bel{D-VM}
d(x,y)\sim \omega_x(\|x-y\|),~~~x,y\in\RN,
\ee
with constants depending only on $n,q,$ and  $\eta$;
\par (b). Let $x,y,z\in\RN$ and let $\lambda\ge 1$ be a constant such that $\|y-z\|\le\lambda\|x-z\|$. Then
\bel{O-1}
d(y,z)\le C\,d(x,z)
\ee
and
\bel{O-2}
\frac{d(x,z)}{\|x-z\|}\le C\,\frac{d(y,z)}{\|y-z\|}
\ee
where $C$ is a constant depending only on $n,q,\eta$, and $\lambda$.
\medskip
\par (c). There exists a constant $C=C(n,q,\eta)>0$ such that for every $x,y\in\RN$ and $z\in(x,y)$ the following inequality
$$
d(x,z)+d(z,y)\le C\, d(x,y)
$$
holds.
\end{claim}
\par {\it Proof.} (a). Let $d=d_q(h)$ where $h\in A_1(\RN)$ and $\|h^q\|_{A_1}\le\eta$. Fix $x\in\RN$. By $v_x$ we denote a function on $\R_+$ such that $v_x(0)=0$ and
$$
v_x(t):=t\left(\frac{1}{|Q(x,t)|}\intl_{Q(x,t)}
h^q\,du\right)^{\frac{1}{q}},~~~~t>0.
$$
Then, by \rf{F-DS},
\bel{D-KL}
d(x,y)\sim v_x(\|x-y\|),~~~y\in\RN.
\ee
Clearly, since $q\ge n$, the function
$$
v_x(t)=2^{-\frac{n}{q}}\,t^{1-\frac{n}{q}}
\left(\,\intl_{Q(x,t)}
h^q\,du\right)^{\frac{1}{q}}
$$
is {\it non-decreasing}. On the other hand, by \rf{4.C}, for every $0<t_1<t_2$
$$
v_x(t_2)/t_2=\left(\frac{1}{|Q(x,t_2)|}\intl_{Q(x,t_2)}
h^q\,du\right)^{\frac{1}{q}}\le \|h^q\|_{A_1} \left(\frac{1}{|Q(x,t_1)|}\intl_{Q(x,t_1)}
h^q\,du\right)^{\frac{1}{q}}
$$
so that
\bel{Q-MN}
v_x(t_2)/t_2\le \eta\,v_x(t_1)/t_1, ~~~~0<t_1<t_2.
\ee
\par It is well known that every non-negative non-decreasing function $v_x$ on $\R_+$ satisfying inequality \rf{Q-MN} is equivalent to a concave non-decreasing function. (In particular, $v_x$ is equivalent to its least concave majorant.) Thus there exists a modulus of continuity $\omega_x\in\MC$ such that $\omega_x(t)\sim v_x(t)$ on $\R_+$ with constants of the equivalence depending only $\eta$.
Combining this equivalence with equivalence \rf{D-KL} we obtain the statement (a) of the claim.\medskip
\par (b). By part (a) of the claim there exists a concave majorant $\omega_z\in\MC$ such that $d(a,z)\sim\omega_z(\|a-z\|)$, $a\in\RN$. Since $\omega_z$ is non-negative concave and non-decreasing, the function $\omega_z(t)/t$ is non-increasing. Hence for every $t_1,t_2>0$ such that $t_1\le \lambda t_2$ we have
$$
\omega_z(t_1)\le \lambda\omega_z(t_2)~~~~\text{and}~~~~
\omega_z(t_2)/t_2\le \lambda\omega_z(t_1)/t_1.
$$
It remains to put $t_1:=\|y-z\|$ and $t_2:=\|x-z\|$ and the part (b) of the claim follows.
\medskip
\par (c). Since $z\in(x,y)$, we have $\|y-z\|,\|x-z\|\le\|x-y\|$ so that by part (b) of the claim
$$
d(x,z)\le C\,d(x,y)~~~\text{and}~~~d(y,z)\le C\,d(x,y)
$$
proving the statement (c) and the claim.\bx
\begin{remark} {\em Equivalence \rf{D-VM} motivates us to refer to the metric space $(\RN,d)$ where $d\in\MQA$ as a {\it variable metric transform} of $\RN$; see Remark \reff{VM-E}. This equivalence shows that given $x\in\RN$ the local behavior of the metric $d$ is similar to that of a certain regular metric transform $d_{\omega_x}:=\omega_x(\|\cdot\|)$ where $\omega_x\in\MC$ is a ``modulus of continuity''. The function $\omega_x$ varies from point to point, and this is the main difference between a regular metric transform (where $\omega_x$ is constant, i.e., the same ``modulus of continuity'' $\omega$ for all $x\in\RN$) and a variable metric transform.
\par Nevertheless, in spite of $\omega_x$ changes, the metric $d\in\MQA$ preserves several important properties of regular metric transforms.
\par For instance, let $E\subset \RN$ be a closed set and let $x\in \RN\setminus E$. Let $\tilde{x}\in E$ be an almost nearest point to $x$ on $E$ with respect to the
Euclidean distance. Then $\tilde{x}$ is an almost nearest to $x$ point with respect to the variable majorant $d$ as well.
\par Another example is the standard Whitney covering of $\RN\setminus E$ by a family of Whitney's cubes. See, e.g. \cite{St}. It is well known that this covering is universal with respect to the family
$$
\Mc\Tc=\{(\RN,d_\omega),~~\omega\in\MC\}
$$
of all metric transforms, i.e., it provides an almost optimal Whitney type extension construction for the family of Lipschitz spaces with respect to metric transforms. As we shall see below the same property holds for variable metric transforms as well.
\par Thus there exists a more or less complete analogy between extension methods for regular and variable metric transforms. In the next sections we present several applications of this approach to extensions of Sobolev functions.\rbx}
\end{remark}
\SECT{3. Sobolev $L^1_p$-space as a union of Lipschitz spaces.}{3}
\addtocontents{toc}{3. Sobolev $L^1_p$-space as a union of Lipschitz spaces. \hfill\thepage\par\VST}
\indent
\par {\it Proof of Theorem \reff{MAIN}.} {\it (Sufficiency).} Let  $n\le q<p$. Let $F\in C(\RN)$ and let $h\in\LPRN$ be a non-negative function such that
$$
|F(x)-F(y)|\le d_q(x,y:h),~~~x,y\in\RN.
$$
Prove that $F\in\LOP$ and
\bel{N-FH}
\|F\|_{\LOP}\le C(n,q,p) \|h\|_{\LPRN}.
\ee
\par Since $d_q\le\delta_q$, for every $x,y\in\RN$ we have
$$
|F(x)-F(y)|\le \delta_q(x,y:h)=\|x-y\| \left(\frac{1}{|Q_{xy}|}\intl_{Q_{xy}} h(u)^q du\right)^{\frac{1}{q}}.
$$
\par Let $Q$ be a cube in $\RN$ and let
$$
F_Q:=\frac{1}{|Q|}\intl_{Q} F(u)\,du.
$$
Then
\be
\frac{1}{|Q|}\intl_{Q} |F(u)-F_Q|\,du
&\le&\sup_{x,y\in Q}|F(x)-F(y)|\le C\sup_{x,y\in Q}
\|x-y\| \left(\frac{1}{|Q_{xy}|}\intl_{Q_{xy}} h(u)^q du\right)^{\frac{1}{q}}\nn\\
&\le& C (\diam Q)^{1-\frac{n}{q}}\sup_{x,y\in Q}
\left(\,\intl_{Q_{xy}} h(u)^q du\right)^{\frac{1}{q}}.\nn
\ee
Clearly, $Q_{xy}\subset 3Q$ for every $x,y\in Q$. Hence
$$
\frac{1}{|Q|}\intl_{Q} |F(u)-F_Q|\,du
\le C(\diam Q)^{1-\frac{n}{q}}
\left(\,\intl_{3Q} h(u)^q du\right)^{\frac{1}{q}}
$$
so that
\bel{QFI}
\frac{1}{|Q|}\intl_{Q} |F(u)-F_Q|\,du\le C\diam Q
\left(\,\frac{1}{|3Q|}\intl_{3Q} h(u)^q du\right)^{\frac{1}{q}}.
\ee
Since $1\le q<p$, for an arbitrary cube $Q\subset\RN$ we have
$$
\frac{1}{|Q|}\intl_{Q} |F(u)-F_Q|\,du\le C\diam Q
\left(\,\frac{1}{|3Q|}\intl_{3Q} h(u)^p du\right)^{\frac{1}{p}}.
$$
\par In \cite{KMc} it is proven that each function $F$ satisfying this condition belongs to $\LOP$ and its gradient $\|\nabla F(x)\|\le C(n,p)h(x)$ a.e. on $\RN$. (See also \cite{H}, p. 266.) This of course implies the required inequality \rf{N-FH}.
\par It is also noted in \cite{H} that for the case $1\le q<p$ inequality \rf{N-FH} can be directly deduced from \rf{QFI} and a theorem of Calder\'{o}n \cite{C1}; see also \cite{CS}. Recall that Calder\'{o}n's theorem states that a function $F\in\LOP$, $1<p<\infty$, provided $F^{\sharp}\in\LPRN$. Here $F^{\sharp}$ is the sharp maximal function of $F$
$$
F^{\sharp}(x):=
\sup_{r>0}\frac1r\,\frac{1}{|Q(x,r)|}\intl_{Q(x,r)} |F(u)-F_{Q(x,r)}|\,du.
$$
Furthermore, $\|F\|_{\LOP}\le C(n,p)\|F^{\sharp}\|_{\LPRN}$.
\par In fact, by \rf{QFI}, for every $x\in\RN$
\be
F^{\sharp}(x)&\le& C\,\sup_{r>0}\frac1r\,\diam Q(x,r)
\left(\,\frac{1}{|3Q(x,r)|}\intl_{3Q(x,r)} h(u)^q du\right)^{\frac{1}{q}}\nn\\
&\le& C\,\sup_{r>0}
\left(\,\frac{1}{|3Q(x,r)|}\intl_{3Q(x,r)} h(u)^q du\right)^{\frac{1}{q}}\le C\left(\Mc[h^q](x)\right)^{\frac{1}{q}}.\nn
\ee
Hence
$$
\|F^{\sharp}\|_{\LPRN}\le C\,\left(\,\,\intl_{\RN}\Mc[h^q]^{\frac{p}{q}}(x)\,dx\right)
^{\frac{1}{p}}.
$$
Since $p/q>1$, by the Hardy-Littlewood maximal theorem,
$$
\|F^{\sharp}\|_{\LPRN}\le C\,\left(\,\,\intl_{\RN}(h^q)^{\frac{p}{q}}(x)\,dx\right)
^{\frac{1}{p}}=C\,\|h\|_{\LPRN}
$$
proving that $F^{\sharp}\in\LPRN$. Hence, by Calder\'{o}n's theorem \cite{C1}, $F\in\LOP$ and
$$
\|F\|_{\LOP}\le C\|F^{\sharp}\|_{\LPRN}\le
C\,\|h\|_{\LPRN}
$$
proving inequality \rf{N-FH} and the sufficiency part of Theorem \reff{MAIN}.
\bigskip
\par {\it (Necessity).} Suppose that $F\in\LOP$. Let $n<q<p$. Fix a constant $\sigma\in(q,p)$, say $\sigma=(p+q)/2$. We put
$h_1(x):=\|\nabla F(x)\|$ and $h_2(x):=\Mc[h_1^\sigma]^{\frac1\sigma}(x),$ $x\in\RN$.  Notice that
$$
h_1=(h_1^\sigma)^{\frac1\sigma}\le \left(\Mc[h_1^\sigma]\right)^{\frac1\sigma}=h_2~~~~
\text{a.e. on}~~ \RN
$$
so that, by \rf{S-P},
$$
|F(x)-F(y)|\leq C(n,q) \|x-y\|\left(\frac{1}{|Q_{xy}|}\intl_{Q_{xy}} h_2^q(u)\,du\right)^{\frac{1}{q}}=
\|x-y\|\left(\frac{1}{|Q_{xy}|}\intl_{Q_{xy}} h^q(u)\,du\right)^{\frac{1}{q}}
$$
where $h=C(n,q)h_2$.
\par Coifman and Rochberg \cite{CR} proved that if
$g\in L_{1,loc}(\RN)$ and $\Mc[g](x)<\infty$ a.e., then  $\Mc[g]^{\theta}\in A_1(\RN)$ for every $0<\theta<1$ and
$\|\Mc[g]^{\theta}\|_{A_1}\le \gamma(n,\theta)$.
\par Let us apply this result to the function $g=h^q=C\Mc[h_1^\sigma]^{\frac{q}{\sigma}}$. Since $0<\theta=q/\sigma<1$, the function $h^q\in A_1(\RN)$ and
$$
\|h^q\|_{A_1}=\|C\Mc[h_1^\sigma]^{\theta}\|_{A_1}
=\|\Mc[h_1^\sigma]^{\theta}\|_{A_1}\le \gamma(n,\theta)=\gamma(n,\tfrac{p+q}{2})=
\eta(n,p,q).
$$
\par Since $p/\sigma>1$, by the Hardy-Littlewood maximal theorem,
\be
\|h\|_{\LPRN}&=& C\left(\,\,\intl_{\RN}\Mc[h_1^\sigma]
^{\frac{p}{\sigma}}(x)\,dx\right)
^{\frac{1}{p}}\le C\left(\,\,\intl_{\RN}(h_1^\sigma)^{\frac{p}{\sigma}}(x)\,dx\right)
^{\frac{1}{p}}\nn\\
&=&C\left(\,\,\intl_{\RN}\|\nabla F\|^p(x)\,dx\right)^{\frac{1}{p}}=C\|F\|_{\LOP}.\nn
\ee
\par Thus  we have proved that for every $q\in(n,p)$ there exists a non-negative function $h\in\LPRN$ such that
\bel{H-C}
h^q\in A_1(\RN),~
\|h^q\|_{A_1}\le\eta(n,p,q), ~\|F\|_{\LOP}\le C\|h\|_{\LPRN}
\ee
and
\bel{F-DL}
|F(x)-F(y)|\le \delta_q(x,y:h),~~~\text{for every}~~x,y\in\RN.
\ee
\par Thus it remains to show the existence of such a function $h$ for $q=n$. Let $h$ be a function satisfying inequalities \rf{H-C} and \rf{F-DL} for $q:=(p+n)/2$. Prove that for some positive constant $C=C(n,p)$ the function $Ch$ satisfies  \rf{H-C} and \rf{F-DL} whenever $q=n$. 
\par In fact, for every cube $Q\subset\RN$
$$
\frac{1}{|Q|}\intl_{Q} h^n(u)\,du\le \left(\frac{1}{|Q|}\intl_{Q} h^q(u)\,du\right)^{\frac{n}{q}}\le
\|h^q\|_{A_1}^{\frac{n}{q}} \left(\essinf\limits_{Q} h^q
\right)^{\frac{n}{q}}=\|h^q\|_{A_1}^{\frac{n}{q}}
\essinf\limits_{Q} h^n
$$
so that $h^n\in A_1(\RN)$ and
$$
\|h^n\|_{A_1}\le\|h^q\|_{A_1}^{\frac{n}{q}}
\le\eta(n,p,q)^{\frac{n}{q}}\le
\eta(n,p,q).
$$
\par Furthermore, by \rf{EW-G},
$$
\dq(x,y:h)\le \|h^q\|_{A_1}^{\frac1q}\,\delta_n(x,y:h)\le C\,\delta_n(x,y:h)
$$
so that, by \rf{F-DL},
$$
|F(x)-F(y)|\le \delta_q(x,y:h)\le C\,\delta_n(x,y:h)=\delta_n(x,y:Ch),~~~~x,y\in\RN.
$$
Since $\|(Ch)^n\|_{A_1}=\|h^n\|_{A_1}$, the function $Ch$ satisfies conditions \rf{H-C} and \rf{F-DL} if $q=n$.
\par Now combining inequality \rf{F-DL} and definition \rf{GD-M} we finally obtain the required estimate \rf{LIP1}.
\par The proof of Theorem \reff{MAIN} is complete.\bx
\SECT{4. Sobolev $L^1_p$-functions on closed subsets of
$\RN$: a proof of Theorem \reff{CR-LOP}.}{4}
\addtocontents{toc}{4. Sobolev $L^1_p$-functions on closed subsets of $\RN$: a proof of Theorem \reff{CR-LOP}. \hfill\thepage\par\VST}
\indent
\par Let $f\in\LOP|_E$ and let $F\in\LOP$ be an arbitrary continuous function such that $F|_E=f$. Given $x\in\RN$ let
\bel{D-SHR}
\fs(x):=\sup_{y,z\in E}\frac{|f(y)-f(z)|}{\|x-y\|+\|x-z\|}.
\ee
Prove that $\fs\in\LPRN$ and
\bel{E-1}
\|\fs\|_{\LPRN}\le C(n,p)\|F\|_{\LOP}.
\ee
\par In fact, let $q:=(n+p)/2$. Given $x\in\RN$ and $y,z\in E$ let $K:=Q(x,\|x-y\|+\|x-z\|)$. Then, by the Sobolev-Poincar\'e inequality \rf{S-P},
\be
|f(y)-f(z)|&=&|F(y)-F(z)|\le C\diam K \left(\frac{1}{|K|} \intl_K \|\nabla F(u)\|^q\,du\right)^{\frac{1}{q}}\nn\\
&\le&
C(\|x-y\|+\|x-z\|)\left(\frac{1}{|K|} \intl_K \|\nabla F(u)\|^q\,du\right)^{\frac{1}{q}}\nn
\ee
where $C=C(n,p)$. Hence
$$
\frac{|f(y)-f(z)|}{\|x-y\|+\|x-z\|}\le C\left(\frac{1}{|K|} \intl_K \|\nabla F(u)\|^q\,du\right)^{\frac{1}{q}}\le C \Mc[\|\nabla F\|^q]^{\frac1q}(x).
$$
Taking the supremum in this inequality over all $y,z\in E$ we obtain
$$
\fs(x)\le C \Mc[\|\nabla F\|^q](x)^{\frac1q},~~~~x\in\RN.
$$
Hence, by the Hardy-Littlewood maximal theorem,
$$
\|\fs\|_{\LPRN}\le C \|\Mc[\|\nabla F\|^q]^{\frac1q}\|_{\LPRN}\le C \|\nabla F\|_{\LPRN}=C\|F\|_{\LOP}
$$
proving \rf{E-1}.
\par Let
$$
I_p(f;E):=\left\{\,\,\intl_{\RN}\left(\,\sup_{y,z\in E}\frac{|f(y)-f(z)|}{\|x-y\|+\|x-z\|}\right)^p\,dx
\right\}^{\frac{1}{p}}.
$$
Then, by inequality \rf{E-1},
$$
I_p(f;E)=\|\fs\|_{\LPRN}\le C\|F\|_{\LOP}.
$$
Taking the infimum in the right hand side of this inequality over all $F\in\LOP$ such that $F|_E=f$ we obtain
$$
I_p(f;E)\le C(n,p)\|f\|_{\LOP|_E}.
$$
\medskip
\par Prove that $f\in \LOP|_E$ and $\|f\|_{\LOP|_E}\le C I_p(f;E)$ provided $f$ is a continuous function on $E$ and $I_p(f;E)<\infty$.
\par Let $\theta:=(q+p)/2$; thus $n<q=(p+n)/2<\theta<p$. Let $h_1:=\Mc[(\fs)^\theta]^{\frac1\theta}$. Clearly, $\fs\le h_1$ a.e. on $\RN$.
\par Let $x,y\in E$ and let $u\in Q_{xy}=Q(x,\|x-y\|)$. Then, by definition \rf{D-SHR},
$$
|f(x)-f(y)|\le\fs(u)(\|u-x\|+\|u-y\|)
$$
so that
$$
|f(x)-f(y)|\le 3\|x-y\|\fs(u)~~~\text{for every}~~u\in Q_{xy}.
$$
Integrating this inequality over the cube $Q_{xy}$ we obtain
$$
|f(x)-f(y)|\le 3\|x-y\| \left(\frac{1}{|Q_{xy}|} \intl_{Q_{xy}} \fs(u)\,du\right)\le 3\|x-y\| \left(\frac{1}{|Q_{xy}|} \intl_{Q_{xy}} h_1(u)\,du\right).
$$
Hence
$$
|f(x)-f(y)|\le 3\|x-y\| \left(\frac{1}{|Q_{xy}|} \intl_{Q_{xy}} h_1^q(u)\,du\right)^{\frac1q}.
$$
so that
\bel{L-A}
|f(x)-f(y)|\le \delta_q(x,y:h_2),~~~~x,y\in E,
\ee
where $h_2:=3h_1$.
\par By a result of Coifman and Rochberg \cite{CR} which we have mentioned in Section 3, the weight
\bel{H-INA}
h_2^q=3^q h_1^q=3^q\Mc[(\fs)^\theta]^{\frac{q}{\theta}}\in A_1(\RN)
\ee
and
\bel{H-A1}
\|h_2^q\|_{A_1}\le \eta(n,p,q).
\ee
Furthermore, since $1<\theta<p$, by the Hardy-Littlewood maximal theorem,
\bel{H-LP}
\|h_2\|_{\LPRN}\le C\|\Mc[(\fs)^\theta]^{\frac{1}{\theta}}\|_{\LPRN}
\le C\|\fs\|_{\LPRN}=C I_p(f;E).
\ee
\par Thus the geodesic distance
$$
d=d_q(x,y:h_2),~~~~x,y\in\RN,
$$
associated with the function $\delta_q(h_2)$, see \rf{GD-M}, belongs to the family $\DPQ$ of metrics defined by \rf{D-LIP}.

\par By Theorem \reff{MAIN-MT},
$$
\delta_q(x,y:h_2)\le C(n,q)\,d_q(x,y:h_2)~~~\text{for all}~~~x,y\in\RN,
$$
so that, by \rf{L-A},
$$
|f(x)-f(y)|\le C\,d_q(x,y:h_2),~~~~x,y\in E.
$$
Hence $f\in\Lip(E;d_q(h_2))$ and $\|f\|_{\Lip(E;d_q(h_2))}\le C$.
\par Let us extend the function $f$ from the set $E$ to all of $\RN$ using the McShane extension formula
$$
F(x)=\inf_{y\in E}\{f(y)+C\,d_q(x,y:h_2)\},~~~~~x\in E.
$$
Then, by McShane's theorem \cite{McS}, the extension $F\in\Lip(\RN;d_q(h_2))$ and
$$
\|F\|_{\Lip(\RN;d_q(h_2))}\le\|f\|_{\Lip(E;d_q(h_2))}\le C.
$$
\par Let $h:=C\,h_2$. We have proved that $F$ satisfies the Lipschitz condition
$$
|F(x)-F(y)|\le C\,d_q(x,y:h_2)=d_q(x,y:h),~~~~x,y\in\RN,
$$
with respect to the metric $d_q(h)$. Furthermore, $h^q=C^q\,h_2^q\in A_1(\RN)$ and  $\|h^q\|_{A_1}=\|h_2^q\|_{A_1}\le\eta(n,p,q)$. Hence, by Theorem \reff{MAIN},
$$
F\in\LOP~~~\text{and}~~~\|F\|_{\LOP}\le C(n,p)\|h\|_{\LPRN}\le C(n,p)\|h_2\|_{\LPRN}
$$
so that, by \rf{H-LP}, $\|F\|_{\LOP}\le C I_p(f;E)$.
\par Since $F\in\LOP$, the function $f=F|_E\in \LOP|_E$ and 
$$
\|f\|_{\LOP|_E}\le \|F\|_{\LOP}\le C I_p(f;E)
$$
proving the theorem.\bx
\begin{remark}\lbl{ALG-MC} {\em Let us present main steps of an extension algorithm which we used in the proof of Theorem \reff{CR-LOP}.
\par Let $q=(n+p)/2$ and $\theta=(q+p)/2$. Let $f$ be a continuous function defined on $E$.\smallskip
\par \textit{Step 1.} We construct the sharp maximal function
$$
\fs(x)= \sup_{y,z\in E}\frac{|f(x)-f(y)|}{\|x-y\|+\|y-z\|},~~~x\in\RN.
$$
\smallskip
\par \textit{Step 2.} We define a weight $h:\RN\to \R_+$ by the formula
$h=\Mc\left[(\fs)^\theta\right]^{\frac1\theta}$.
\smallskip
\par \textit{Step 3.} Using formula \rf{GD-M} we construct the geodesic distance $d=d_q(x,y:h)$ associated with the ``pre-metric''
$$
\delta_q(x,y:h)= \|x-y\|\left(\frac{1}{|Q_{xy}|}
\intl_{Q_{xy}}h^q(u)du\right)^{\frac{1}{q}},~~~~x,y\in\RN.
$$
\smallskip
\par \textit{Step 4.} We extend the function $f:E\rightarrow \R$ to a function $F$ defined on all of $\RN$ using the McShane's extension formula
$$
F(x)=\inf_{y\in E}\left\{f(y)+C\,d_q(x,y:h)\right\}.
$$
where $C$ is the constant from inequality \rf{W-G}.
\par By Theorem \reff{CR-LOP}, the function $F$ provides an almost optimal extension of the function $f$ to a function from the Sobolev space $\LOP$ whenever $\fs\in\LPRN$.
\par Of course, the extension operator $f\mapsto F$ is non-linear. Notice that this algorithm is new even for families of ``nice'' sets $E\subset\RN$ (like, say,  closures of Lipschitz domains etc.).
\par As we have mentioned above, the classical Whitney extension method also provides an almost optimal extensions of functions from the trace space $\LOP|_E$ to functions from $\LOP$. See \cite{Sh2}. We give an alternative proof of this property of the Whitney extension operator in the next section.\rbx}
\end{remark}
\SECT{5. Sobolev-Poincar\'e inequality and extensions of
$L^m_p$-functions.}{5}
\addtocontents{toc}{5. Sobolev-Poincar\'e inequality and extensions of $L^m_p$-functions. \hfill\thepage\par\VST}
\indent
\par In this section we generalize the approach presented in the previous sections to the case of the Sobolev spaces $\LMP$ whenever $p>n$. In particular, we shall prove  representation \rf{DM-PR}.
\par First let us recall the \SP inequality for $\LMP$-functions whenever $p>n$. Let $q\in(n,p)$ and let $F\in \LMP$. Given $x\in \RN$ let $P_x:=T^{m-1}_x[F]$ be the Taylor polynomial of $F$ at $x$ of order $m-1$. Then for every  $x,y\in \RN$ and every multiindex $\beta$, $|\beta|\le m-1,$
\bel{SP-M}
|D^{\beta}P_x(x)-D^{\beta}P_y(x)|\le C\,\|x-y\|^{m-|\beta|}\left(\frac{1}{|Q_{xy}|}
\intl_{Q_{xy}}(\nabla ^mF(u))^qdu\right)^{\frac{1}{q}}.
\ee
In particular, for every $\alpha$, $|\alpha|=m-1,$ the following inequality
$$
|D^{\alpha}F(x)-D^{\alpha}F(y)|\le C\,\|x-y\|\left(\frac{1}{|Q_{xy}|}\intl_{Q_{xy}}(\nabla ^mF(u))^qdu\right)^{\frac{1}{q}},~~~x,y\in\RN,
$$
holds. Here $C=C(n,m,q)$. Recall also that $Q_{xy}=Q(x,\|x-y\|)$, and $\nabla^mF$ is defined by \rf{N-M}.
\medskip
\par Let $d$ be a metric on $\RN$. In Section 1 we have introduced the space $\CMD$ which consists of all $C^m-$functions $F$ on $\RN$ satisfying the following condition: there exists a constant $\lambda=\lambda(F)>0$ such that for every $\alpha$, $|\alpha|=m$, and every $x,y\in\RN$
$$
|D^{\alpha}F(x)-D^{\alpha}F(y)|\le \lambda\, d(x,y).
$$
\par In other words, $F\in \CMD$ if its partial derivatives $D^{\alpha}F$ of order $|\alpha|=m$ belong to $\Lip(\RN;d)$. The space $\CMD$ is normed by
\bel{D-CD}
\|F\|_{\CMD}=\sbig_{|\alpha|=m}\,\,\sup_{x,y\in \RN, x\ne y} \frac{|D^{\alpha}F(x)-D^{\alpha}F(y)|}{d(x,y)}.
\ee
\par We say that a metric $d$ on $\RN$ is \textit{pseudoconvex} if there exists a constant $\lambda_d\ge 1$ such that for every $x,y,z\in \RN$, $z\in (x,y)$, the following inequality
\bel{B-3}
d(x,z)+d(x,y)\le \lambda_d\,d(x,y)
\ee
holds.
\par Let $F\in C^m(\RN)$, $y\in\RN$, and let
$$
T^m_y[F](x)=\smed_{|\alpha|=m}\,\,\frac{1}{\alpha !}\, D^{\alpha}
F(y)(x-y)^{\alpha},~~~~x\in\RN,
$$
be the Taylor polynomial of $F$ of degree $m$ at $y$. The following claim is a variant of the Taylor formula for the space $\CMD$ whenever $d$ is a pseudoconvex metric. 
\begin{claim}\lbl{Cl7.1} Let $d$ be a pseudoconvex metric on $\RN$. Then for every $F\in \CMD$ and every multi-index $\beta$ with $|\beta|\le m$ the following inequality
$$
|D^{\beta}F(x)-D^{\beta}\left(T^m_y[F]\right)(x)|\le C\|F\|_{\CMD}\|x-y\|^{m-|\beta|}d(x,y),~~~~x,y\in \RN,
$$
holds. Here $C=C(n,m,\lambda_d)$.
\end{claim}
\par {\it Proof.} For $|\beta|=m$ the claim follows from definition \rf{D-CD} of the norm in the space $\CMD$, so we can assume that $|\beta|<m$.
\par  Let us make use of the following well known identity:
$$
F(x)=T^m_y[F](x)+m\,\sbig_{|\alpha|=m}\,\frac{1}{\alpha !}(x-y)^{\alpha}\int_0^1(1-t)^{m-1}(D^{\alpha}F(x+t(x-y))-
D^{\alpha}F(y))dt
$$
provided $m>0$, $x,y\in \RN$, and $F\in C^m(\RN)$.
\par Let us apply this identity to $D^{\beta}F$. We obtain
\be
D^{\beta}F(x)&=&
T^{m-|\beta|}_y[D^{\beta}F](x)
+(m-|\beta|)\smed_{|\alpha|=m-|\beta|}\,
\frac{1}{\alpha!}(x-y)^{\alpha}\nn\\
&\cdot&
\int^1_0(1-t)^{m-|\beta|-1}
\left\{D^{\alpha}(D^{\beta}F)(x+t(x-y))
-D^{\alpha}(D^{\beta}F)(y)\right\}dt\,.\nn
\ee
Since
$$
T^{m-|\beta|}_y[D^{\beta}F](x)=D^{\beta}
\left(T^m_y[F]\right)(x)~~~\text{for every}~~~x,y\in \RN,
$$
we have
\be
&&D^{\beta}F(x)-D^{\beta}(T^m_y[F])(x)
=(m-|\beta|)\smed_{|\alpha|=m-|\beta|}
\frac{1}{\alpha!}\,(x-y)^{\alpha}\nn\\
&\cdot&\int^1_0(1-t)^{m-|\beta|-1}
\left\{D^{\alpha+\beta}F)(x+t(x-y))
-D^{\alpha+\beta}F)(y)\right\}dt.\nn
\ee
Hence
$$
|D^{\beta}F(x)-D^{\beta}(T^m_y[F])(x)|\le
m\smed_{|\alpha|+|\beta|=m}
\,\|x-y\|^{|\alpha|}
\sup_{z\in [x,y]}|D^{\alpha+\beta}F)(z)
-D^{\alpha+\beta}F)(y)|.
$$
Combining this inequality with definition \rf{D-CD} we obtain
$$
|D^{\beta}F(x)-D^{\beta}(T^m_y[F])(x)|\le
m\|x-y\|^{m-|\beta|}\|F\|_{\CMD}
\sup_{z\in [x,y]}d(z,y).
$$
Since $d$ is pseudoconvex, by \rf{B-3},
$$
|D^{\beta}F(x)-D^{\beta}(T^m_y[F])(x)|\le
\lambda_d\,m\, \|x-y\|^{m-|\beta|}\,\|F\|_{\CMD}\,d(x,y)
$$
proving the claim. \bx\medskip
\par Let us apply this result  to a metric which belongs to the family $\DPQ$ whenever $q\in[n,p)$, see \rf{D-LIP}. Notice that, by part (c) of Claim \reff{VMT}, every metric $d\in\DPQ$ is quasiconvex. This property of $d$ and Claim \reff{Cl7.1} imply the following
\begin{proposition}\lbl{Prop-FT} Let $d\in\DPQ$ where $q\in[n,p)$. Then for every $F\in \CMD$, every $x,y\in \RN$ and every $\beta$, $|\beta|\le m$, the following inequality
$$
|D^{\beta}F(x)-D^{\beta}(T^m_y[F])(x)|\le
C\|F\|_{\CMD}\|x-y\|^{m-|\beta|}d(x,y)
$$
holds. Here $C=C(n,q,p)$.
\end{proposition}
\medskip\par Our next result is a generalization of Theorem \reff{MAIN} to the case of the Sobolev space $\LMP$.
\begin{theorem}\lbl{SM-MAIN} Let $m$ be a positive integer and let $n\le q<p<\infty$. There exists a positive constant $\eta=\eta(n,p,q)$ depending only on $n,q,$ and $p$ such that the following statement is true:
\par A $C^{m-1}$-function $F$ belongs to $\LMP$ if and only if there exists a non-negative function $h\in\LPRN$ such that $h^q\in \ARN$, $\|h^q\|_{A_1}\le \eta,$ and for every multiindex $\beta$, $|\beta|=m-1$, and every $x,y\in\RN$ the following inequality
\bel{LIP2}
|D^\beta F(x)-D^\beta F(y)|\le d_q(x,y:h)
\ee
holds. Furthermore,
$$
\|F\|_{\LMP}\sim \inf \|h\|_{\LPRN}.
$$
The constants of this equivalence depend only on $m,n,q,$ and $p$.
\end{theorem}
\par {\it Proof.} The case $m=1$ is proven in Theorem \reff{MAIN}. The case $m>1$ easily follows from this result.
\par {\it (Necessity.)} Suppose that $F\in\LMP$. Then $D^\beta F\in\LOP$ for every $\beta$ with $|\beta|=m-1$ and 
\bel{DBF}
\|D^\beta F\|_{\LOP}\le C\,\|F\|_{\LMP}.
\ee
\par By Theorem \reff{MAIN}, there exists a function $h_\beta\in \LPRN$ satisfying the following conditions:
$h_\beta^q\in \ARN$, $\|h_\beta^q\|_{A_1}\le \eta=\eta(n,p,q),$
\bel{LP-B}
\|h_\beta\|_{\LPRN}\le C\,\|D^\beta F\|_{\LOP},
\ee
and
\bel{B-W}
|D^\beta F(x)-D^\beta F(y)|\le d_q(x,y:h_\beta),~~~~x,y\in\RN.
\ee
\par Let
$$
h:=\left(\,\smed_{|\beta|=1}\,\,h^q_\beta\right)^{\frac{1}{q}}.
$$
\par Since $h_\beta\le h$, by definitions \rf{DQ} and \rf{GD-M}, $d_q(h_\beta)\le d_q(h)$ for every $\beta, |\beta|=m-1$. This inequality and \rf{B-W} imply \rf{LIP2}. On the other hand, by \rf{DBF} and \rf{LP-B}, $\|h\|_{\LPRN}\le C\|F\|_{\LMP}$.
\par It remains to prove that $h^q\in \ARN$. In fact, since $\|h^q_\beta\|_{A_1}\le \eta,$ $|\beta|=m-1,$ by definition \rf{4.3}, for every cube $Q\subset\RN$ we have  
\be
\frac{1}{|Q|}\int_Q h^q(u)\,du&=&\smed_{|\beta|=m-1}  \frac{1}{|Q|}\int_Q h_\beta^q(u)\,du\le \eta\smed_{|\beta|=m-1} \essinf\limits_Q h_\beta^q\nn\\
&\le&\eta\essinf\limits_Q\left(\smed_{|\beta|=m-1}
h_\beta^q\right)=\eta\essinf\limits_Q h^q.\nn
\ee
Hence $\|h\|_{A_1}\le \eta$ proving the necessity.\medskip
\par {\it (Sufficiency.)} Suppose that $F\in \CMON$ and there exists a function $h\in\LPRN$ such that
$h^q\in \ARN$, $\|h^q\|_{A_1}\le \eta,$ and for every $\beta, |\beta|=m-1,$ inequality \rf{LIP2} is satisfied.  Then, by Theorem \reff{MAIN}, $D^\beta F\in\LOP$  and
$$
\|D^\beta F\|_{\LOP}\le C\|h\|_{\LPRN}
$$
so that $F\in \LMP$ and $\|F\|_{\LMP}\le C\|h\|_{\LPRN}$.
\par The proof of Theorem \reff{SM-MAIN} is complete. \bx\medskip
\begin{remark}\lbl{CM-L} {\em It is obvious that Theorem \reff{SM-MAIN} implies representation \rf{DM-PR}. Furthermore,  by this theorem, for each $F\in\LMOP$
$$
\|F\|_{\LMOP}\sim \inf\{\|h\|_{\LPRN}: d=d_q(h)\in\DPQ, \|F\|_{\CMD}\le 1\}
$$
with constants in this equivalence depending only on $n,q,$ and $p$.\rbx}
\end{remark}
\bigskip
\par {\bf 5.1 A Whitney-type extension theorem for the space $\CMD$.}
\addtocontents{toc}{~~~~5.1 A Whitney-type extension theorem for the space $\CMD$. \hfill \thepage\par}
\medskip
\par We turn to the proof of Theorem \reff{JET-S}. This proof is based on representation \rf{DM-PR} and the following Whitney-type extension theorem for spaces $\CMD$ where $d\in\DPQ$ is a metric transform.
\begin{theorem} \lbl{CMD-E} Let $p\in(n,\infty)$, $q\in[n,p)$, $m\in\N$, and let $d\in\DPQ$. Let
$$
J=\{P_x\in\PMP: x\in E\}
$$
be a polynomial field defined on a closed set $E\subset\RN$. There exists a function $F\in\CMD$ such that $T^{m}_x[F]=P_x$ for every $x\in E$ if and only if the following quantity
$$
\Lc_{m,d}(J):=\sbig_{|\alpha|\le m}\,\,\,
\sup_{x,y\in E,\,x\ne y}\,\, \frac{|D^{\alpha}P_x(x)-D^{\alpha}P_y(x)|}
{\|x-y\|^{m-|\alpha|}\,d(x,y)}
$$
is finite. Furthermore,
$$
\Lc_{m,d}(J)\sim\inf\left\{\|F\|_{\CMD}:F\in \CMD,\,\, T^{m}_x[F]=P_x~~\text{for all}~~x\in E\right\}
$$
with constants depending only on $n,p,q,$ and $m$.
\end{theorem}
\par {\it Proof.} Suppose that the quantity
$$
\Ic(J):=\inf\left\{\|F\|_{\CMD}:F\in \CMD,\,\, T^{m}_x[F]=P_x~~\text{for every}~~x\in E\right\}
$$
is finite. Then, by Proposition \reff{Prop-FT}, for every function $F\in \CMD$ such that $T^{m}_x[F]=P_x$, $x\in E$, the following inequality
$$
\Lc_{m,d}(J)\le C\, \|F\|_{\CMD}
$$
holds. Taking the infimum in this inequality over all such function $F$ we obtain that $\Lc_{m,d}(J)\le C\, \Ic(J)$ with $C=C(n,p,q,m)$.
\medskip
\par Prove that  $\Ic(J)\le C\, \Lc_{m,d}(J)$ where $C$ is a constant depending only on $n,p,q,$ and $m$.
\par Let $\lambda>0$ and let $J=\{P_x\in\PMP: x\in E\}$ be a polynomial field on $E$ such that for every multiindex $\alpha$, $|\alpha|\le m$, and every $x,y\in E$ we have
\bel{W-P}
|D^{\alpha}P_x(x)-D^{\alpha}P_y(x)|\le \lambda\,\|x-y\|^{m-|\alpha|}\,d(x,y).
\ee
\par Let $F:\RN\to\R$ be an extension of the field $J$ obtained by the Whitney extension method. Prove that $F$ belongs to the space $\CMD$ and the norm of $F$ in this space does not exceed $C(n,q,p,m)\lambda$.
\par But first let us recall the Whitney extension construction. Since $E$ is a closed set, the set $\RN\setminus E$ is open so that it admits a Whitney covering $W_E$ by a family $W_E$ of non-overlapping cubes. See, e.g. \cite{St}, or \cite{Guz}. These cubes have the following properties:
\medskip
\par (i). $\RN\setminus E=\cup\{Q:Q\in W_E\}$;\smallskip
\par (ii). For every cube $Q\in W_E$ we have
\bel{DQ-E}
\diam Q\le \dist(Q,E)\le 4\diam Q.
\ee
\par We are also needed certain additional properties of
Whitney's cubes which we present in the next lemma. These
properties easily follow from constructions of the Whitney covering given in \cite{St} and \cite{Guz}.
\par Given a cube $Q\subset\RN$ let $Q^*:=\frac{9}{8}Q$.
\begin{lemma}\lbl{Wadd}
(1). If $Q,K\in W_E$ and $Q^*\cap K^*\ne\emptyset$, then
$$
\frac{1}{4}\diam Q\le \diam K\le 4\diam Q\,;
$$
\par (2). For every cube $K\in W_E$ there are at most
$N=N(n)$ cubes from the family
$W_E^*:=\{Q^*:Q\in W_E\}$
which intersect $K^*$;
\medskip
\par (3). If $Q,K\in W_E$, then $Q^*\cap K^*\ne\emptyset$
if and only if  $Q\cap K\ne\emptyset$.
\end{lemma}
\par Let $\Phi_E:=\{\varphi_Q:Q\in W_E\}$ be a smooth partition of unity subordinated to the Whitney decomposition $W_E$. Recall the main properties of this partition.
\begin{lemma}\lbl{P-U} The family of functions $\Phi_E$ has the following properties:
\medskip
\par (a). $\varphi_Q\in C^\infty(\RN)$ and
$0\le\varphi_Q\le 1$ for every $Q\in W_E$;\smallskip
\medskip
\par (b). $\supp \varphi_Q\subset Q^*(:=\frac{9}{8}Q),$
$Q\in W_E$;\smallskip
\medskip
\par (c). $\smed\,\{\varphi_Q(x):Q\in W_E\}=1$ for every
$x\in\RN\setminus S$;\smallskip
\medskip
\par (d). For every cube $Q\in W_E$, every $x\in\RN$ and every multiindex $\beta, |\beta|\le m,$ the following inequality
$$
|D^\beta\varphi_Q(x)| \le C(n,m)\,(\diam Q)^{-|\beta|}
$$
holds.
\end{lemma}
\medskip
\par Given a cube $Q\in W_E$ by $a_Q$ we denote a point nearest to $Q$ on the set $E$. Notice the following property of $a_Q$ which follows from inequality \rf{DQ-E}:
\bel{AQ}
a_Q\in 9 Q~~~\text{for every cube}~~~Q\in W_E.
\ee
\par By $P^{(Q)}$ we denote the polynomial $P_{a_Q}$. Finally, we define the extension $F$ by the Whitney extension formula:
\bel{DEF-F}
F(x):=\left \{
\begin{array}{ll}
P_x(x),& x\in E,\\\\
\sbig\limits_{Q\in W_E}\,\,
\varphi_Q(x)P^{(Q)}(x),& x\in\RN\setminus E.
\end{array}
\right.
\ee
\medskip
\par Let us note that the metric $d$ is continuous with respect to the Euclidean distance, i.e., for every $x\in\RN$
\bel{C-M}
d(x,y)\to 0 ~~~\text{as}~~~\|x-y\|\to 0.
\ee
In fact, since $d\in\DPQ$, by definition \rf{D-LIP}, $d=d_q(h)$ where $h\in L_{q,loc}(\RN)$ is a non-negative function such that $h^q\in A_1(\RN)$ and $\|h^q\|_{A_1}\le\eta(n,p,q)$. Then, by Theorem \reff{MAIN-MT},
$$
d(x,y)\sim\dq(x,y:h)=\|x-y\| \left(\frac{1}{|Q_{xy}|}\intl_{Q_{xy}} h(u)^q du\right)^{\frac{1}{q}}=2^{\frac{n}{q}}
\|x-y\|^{1-\frac{n}{q}}
\left(\intl_{Q_{xy}} h(u)^q du\right)^{\frac{1}{q}}.
$$
Since $n\le q$, $h\in L_{q,loc}(\RN)$ and $\diam Q_{xy}=2\,\|x-y\|\to 0$ as $\|x-y\|\to 0$, we have $d(x,y)\to 0$ proving \rf{C-M}.\medskip
\par Hence, by \rf{W-P}, for every multiindex $\beta$, $|\beta|\le m$, we have
$$
D^{\beta}P_x(x)-D^{\beta}P_y(x)
=o( \|x-y\|^{m-|\beta|}),~~~~~x,y\in E,
$$
so that the $m$-jet $\{P_x\in\PMP: x\in E\}$ satisfies the hypothesis of the Whitney extension theorem \cite{W1}. By this theorem the extension $F:\RN\to\R$ defined by the formula \rf{DEF-F} is a $C^m$-function such that
\bel{W-EF}
D^{\beta}F(x)=D^{\beta}P_x(x)~~\text{for every}~~x\in E~~\text{and every}~~\beta, |\beta|\le m.
\ee
\par Prove that $\|F\|_{\CMD}\le C\lambda$, i.e., for every multiindex $\alpha$, $|\alpha|=m$, and every $x,y\in \RN$ the following inequality
\bel{AIM}
|D^{\alpha}F(x)-D^{\alpha}F(y)|\le C\lambda\,d(x,y).
\ee
holds.
\par Consider four cases.
\par {\it The first case: $x,y\in E$.} Since $P_y\in\PMP$, for every multiindex $\alpha$ of order $m$ the function $D^\alpha P_y$ is a constant function. In particular, $D^{\alpha}P_y(x)=D^{\alpha}P_y(y)$. Hence, by \rf{W-EF},
$$
|D^{\alpha}F(x)-D^{\alpha}F(y)|=
|D^{\alpha}P_x(x)-D^{\alpha}P_y(y)|
=|D^{\alpha}P_x(x)-D^{\alpha}P_y(x)|
$$
so that, by \rf{W-P},
$$
|D^{\alpha}F(x)-D^{\alpha}F(y)|\le \lambda\,d(x,y)
$$
proving \rf{AIM} in the case under consideration.\bigskip
\par {\it The second case: $x\in E$, $y\in\RN\setminus E$.} Given a Whitney cube $K\in W_E$ let
\bel{TK}
T(K):=\{Q\in W_E: Q\cap K\ne\emp\}
\ee
be a family of all Whitney cubes touching $K$.
\begin{lemma}\lbl{C-A} Let $K\in W_E$ be a Whitney cube and let $y\in K^*=\tfrac98 K.$ Then for every multiindex $\alpha$ the following inequality
$$
|D^{\alpha}F(y)-D^{\alpha}P_{a_K}(y)|\le C\,
\smed\limits_{Q\in T(K)}\,\,\smed\limits_{|\xi|\le m}\,\,
(\diam K)^{|\xi|-|\alpha|} \,|D^{\xi}P_{a_Q}(a_K)-D^{\xi}P_{a_K}(a_K)|
$$
holds. Here $C=C(n,m,\alpha)$.
\end{lemma}
\par {\it Proof.} We notice that, by part (2) of Lemma \reff{Wadd}, $\#T(K)\le N(n)$, and, by part (3) of this lemma,
$$
T(K)=\{Q\in W_E: Q^*\cap K^*\ne\emp\}.
$$
\par Recall $P^{(Q)}=P_{a_Q}$ and $a_Q\in 9Q$ for every $Q\in W_E$. Let us estimate the quantity
$$
I:=|D^{\alpha}F(y)-D^{\alpha}P_{a_K}(y)|.
$$
By formula \rf{DEF-F} and by part (c) of Lemma \reff{Wadd}, 
$$
F(y)-P_{a_K}(y)=\smed\limits_{Q\in W_E}
\varphi_Q(y)(P^{(Q)}(y)-P_{a_K}(y))
$$
so that, by part (b) of Lemma \reff{Wadd} and by definition \rf{TK},
$$
F(y)-P_{a_K}(y)=\smed\limits_{Q\in T(K)}
\varphi_Q(y)(P^{(Q)}(y)-P_{a_K}(y)).
$$
Hence
$$
I:=|D^\alpha F(y)-D^\alpha P_{a_K}(y)|\le \smed\limits_{Q\in T(K)}
|D^\alpha(\varphi_Q(y)(P^{(Q)}(y)-P_{a_K}(y)))|
$$
so that
\bel{I2}
I\le \smed\limits_{Q\in T(K)} A_Q(y;\alpha)
\ee
where
$$
A_Q(y;\alpha):=
|D^\alpha(\varphi_Q(y)(P^{(Q)}(y)-P_{a_K}(y)))|.
$$
\par Let $Q\in T(K)$. Then
$$
A_Q(y;\alpha)\le C\smed\limits_{|\beta|+|\gamma|=|\alpha|}
|D^\beta\varphi_Q(y)|\,|D^\gamma(P^{(Q)}(y)-P_{a_K}(y))|
$$
so that, by part(d) of Lemma \reff{P-U},
$$
A_Q(y;\alpha)\le C\smed\limits_{|\beta|+|\gamma|=|\alpha|}
(\diam Q)^{-|\beta|} \,|D^\gamma(P^{(Q)}(y)-P_{a_K}(y))|.
$$
Since $Q\cap K\ne\emp$, by part(1) of Lemma \reff{Wadd}, $\diam Q\sim\diam K$ so that
\bel{AQA}
A_Q(y;\alpha)\le C\smed\limits_{|\beta|+|\gamma|=|\alpha|}
(\diam K)^{-|\beta|} \,|D^\gamma(P^{(Q)}(y)-P_{a_K}(y))|.
\ee
\par Let us estimate the distance between $a_Q$ and $y$. By \rf{DQ-E},
$$
\|a_Q-y\|\le \diam K^*+\diam Q+\dist(Q,E)
\le 2\diam K+\diam Q+4\diam Q.
$$
Since $Q\cap K\ne\emp$, by part (1) of Lemma \reff{Wadd},
\bel{D-AQ}
\|a_Q-y\|\le 2\diam K+4\diam K+16\diam K=22\diam K.
\ee
\par Let
$$
\tP_Q:=P^{(Q)}-P_{a_K}=P_{a_Q}-P_{a_K}.
$$
Let us estimate the quantity $|D^\gamma\tP_Q(y)|$. Since $\tP_Q\in\PMP$, we can represent this polynomial as
$$
\tP_Q(z)=\smed\limits_{|\xi|\le m}\,\,\frac{1}{\xi!}\, D^{\xi}\tP_Q(a_K)\,(z-a_K)^{\xi}.
$$
Hence
$$
D^\gamma\tP_Q(z)=\smed\limits_{|\gamma|\le |\xi|\le m}\,\,\frac{1}{(\xi-\gamma)!}\, D^{\xi}\tP_Q(a_K)\,(z-a_K)^{\xi-\gamma}
$$
so that
$$
|D^\gamma\tP_Q(y)|\le C\smed\limits_{|\gamma|\le |\xi|\le m} \,
|D^{\xi}\tP_Q(a_K)|\,\|y-a_K\|^{|\xi|-|\gamma|}\le
C\smed\limits_{|\gamma|\le |\xi|\le m} \,\,
(\diam K)^{|\xi|-|\gamma|}\,|D^{\xi}\tP_Q(a_K)|.
$$
Combining this inequality with \rf{AQA} we obtain
\be
A_Q(y;\alpha)&\le& C\smed\limits_{|\beta|+|\gamma|=|\alpha|}
(\diam K)^{-|\beta|} \,|D^\gamma\tP_Q(y)|\nn\\
&\le& C\smed\limits_{|\beta|+|\gamma|=|\alpha|}
(\diam K)^{-|\beta|} \,\smed\limits_{|\gamma|\le |\xi|\le m}(\diam K)^{|\xi|-|\gamma|}\,|D^{\xi}\tP_Q(a_K)|\nn\\
&\le&
C \,\smed\limits_{|\xi|\le m}(\diam K)^{|\xi|-|\alpha|}\,|D^{\xi}\tP_Q(a_K)|.\nn
\ee
Hence, by \rf{I2},
$$
I\le C\,\smed\limits_{Q\in T(K)}\,\,
\smed\limits_{|\xi|\le m}
(\diam K)^{|\xi|-|\alpha|}\,|D^{\xi}\tP_Q(a_K)|
$$
proving the lemma.\bx
\begin{lemma}\lbl{C-2} Let $x\in E$ and let $K\in W_E$ be a Whitney cube. Then for every $y\in K$ and every $\alpha,|\alpha|=m,$ the following inequality
\be
|D^{\alpha}F(x)-D^{\alpha}F(y)|&\le& C(n,m)\,\{
|D^{\alpha}P_x(x)-D^{\alpha}P_{a_K}(x)|\nn\\
&+&
\smed\limits_{Q\in T(K)}\,\,\smed\limits_{|\xi|\le m}
(\diam K)^{|\xi|-m} \,|D^{\xi}P_{a_Q}(a_K)-D^{\xi}P_{a_K}(a_K)|\}\nn
\ee
holds.
\end{lemma}
\par {\it Proof.} We have
$$
|D^{\alpha}F(x)-D^{\alpha}F(y)|
\le |D^{\alpha}F(x)-D^{\alpha}P_{a_K}(y)|
+|D^{\alpha}P_{a_K}(y)-D^{\alpha}F(y)|=I_1+I_2.
$$
Since $P_{a_K}\in\PMP$ and $|\alpha|=m$, the function $D^\alpha P_y$ is a constant function so that
$$
D^{\alpha}P_{a_K}(y)=D^{\alpha}P_{a_K}(x).
$$
Since $x\in E$, by \rf{W-EF}, $D^{\alpha}F(x)=D^{\alpha}P_{x}(x)$ so that, by \rf{W-P},  
$$
I_1:=|D^{\alpha}F(x)-D^{\alpha}P_{a_K}(y)|=
|D^{\alpha}P_{x}(x)-D^{\alpha}P_{a_K}(x)|.
$$
It remains to apply Lemma \reff{C-A} to the quantity
$I_2:=|D^{\alpha}F(y)-D^{\alpha}P_{a_K}(y)|$,
and the lemma follows.\bx
\bigskip
\par Let us prove inequality \rf{AIM} for arbitrary $x\in E$, $y\in\RN\setminus E$ and $\alpha,\,|\alpha|=m$. Let $y\in K$ for some $K\in W_E$. By Lemma \reff{C-2},
$$
|D^{\alpha}F(x)-D^{\alpha}F(y)|\le C\{J_1+J_2\}
$$
where
$$
J_1:=|D^{\alpha}P_x(x)-D^{\alpha}P_{a_K}(x)|
$$
and
$$
J_2:=
\smed\limits_{Q\in T(K)}\,\,\smed\limits_{|\xi|\le m}
(\diam K)^{|\xi|-m} \,|D^{\xi}P_{a_Q}(a_K)-D^{\xi}P_{a_K}(a_K)|
$$
\par First let us estimate $J_2$. By \rf{W-P}, for every $\xi, |\xi|\le m,$ and every $Q\in T(K)$, we have
$$
|D^{\xi}P_{a_Q}(a_K)-D^{\xi}P_{a_K}(a_K)|
\le
\lambda\,\|a_Q-a_K\|^{m-|\xi|}d(a_Q,a_K).
$$
Hence
$$
J_2\le\lambda
\smed\limits_{Q\in T(K)}\,\,\smed\limits_{|\xi|\le m}
(\diam K)^{|\xi|-m}\|a_Q-a_K\|^{m-|\xi|}d(a_Q,a_K).
$$
But, by \rf{D-AQ},
\bel{A-4}
\|a_Q-a_K\|\le \|a_Q-y\|+\|y-a_K\|\le 23\diam K
\ee
proving that
\bel{E-J2}
J_2\le C(n,m)\,\lambda
\smed\limits_{Q\in T(K)}d(a_Q,a_K).
\ee
\medskip
\par Now prove that for some constant $C=C(n,p,q)$
\bel{D-AY}
d(a_Q,y)\le C\,d(x,y)~~~\text{for every}~~~Q\in T(K).
\ee
In fact, by \rf{D-AQ}, $\|a_Q-y\|\le 22\diam K$. Since $x\in E$ and $y\in K$, by \rf{DQ-E},
$$
\diam K\le 4\dist(K,E)\le 4\|x-y\|
$$
so that $\|a_Q-y\|\le 88\|x-y\|$. Hence, by part (b) of Claim \reff{VMT}, see inequality \rf{O-1},
$d(a_Q,y)\le C(n,p,q)\,d(x,y)$ proving \rf{D-AY}.
\medskip
\par Now we have
$$
d(a_Q,a_K)\le d(a_Q,y)+d(y,a_K)\le C\,d(x,y)
$$
so that, by \rf{E-J2},
$$
J_2\le C\,\lambda\, \#T(K)\,d(x,y)\le C\,\lambda\, d(x,y).
$$
See part (2) of Lemma \reff{Wadd}.
\par On the other hand, by \rf{W-P} and \rf{D-AY},
$$
J_1:=|D^{\alpha}P_x(x)-D^{\alpha}P_{a_K}(x)|\le \lambda \,d(x,a_K)\le \lambda (d(x,y)+d(y,a_K))\le C\,\lambda\,d(x,y).
$$
Finally,
$$
|D^{\alpha}F(x)-D^{\alpha}F(y)|\le C\,\{J_1+J_2\}\le C\lambda\,d(x,y).
$$
\bigskip
\par {\it The third case: $y\in K$, $K\in W_E$ and $x\in\RN\setminus K^*$.}
\par Since $K^*=\tfrac98 K$ and $x\notin K^*$, we have
$$
\|x-y\|\ge \tfrac{1}{16}\diam K.
$$
Let $a\in E$ be a point nearest to $x$ on $E$. Then
$$
\|a-x\|=\dist(x,E)\le\dist(y,E)+\|x-y\|\le \dist(K,E)+
\diam K+\|x-y\|
$$
so that, by \rf{DQ-E},
$$
\|a-x\|\le 4\diam K+\diam K+\|x-y\|\le 81\|x-y\|.
$$
Hence, by part (b) of Claim \reff{VMT}, see \rf{O-1}, $d(a,x)\le C\,d(x,y)$.
\par We have
$$
\|y-a\|\le \|x-y\|+\|x-a\|\le 82\|x-y\|
$$
so that again, by \rf{O-1}, $d(y,a)\le C\,d(x,y)$. We obtain
$$
|D^{\alpha}F(x)-D^{\alpha}F(y)|\le |D^{\alpha}F(x)-D^{\alpha}F(a)|+
|D^{\alpha}F(a)-D^{\alpha}F(y)|
$$
so that, by the result proven in the second case,
$$
|D^{\alpha}F(x)-D^{\alpha}F(y)|\le C\,\lambda\,(d(x,a)+d(y,a))\le C\,\lambda\,d(x,y).
$$
\bigskip
\par {\it The fourth case: $y\in K$, $x\in K^*$ where $K\in W_E$.} The proof of inequality \rf{AIM} in this case is based on the next
\begin{lemma}\lbl{C-4} Let $K\in W_E$ be a Whitney cube and let $x,y\in K^*$. Then for every multiindex  $\alpha,|\alpha|=m,$ the following inequality
$$
|D^{\alpha}F(x)-D^{\alpha}F(y)|\le C\,
\frac{\|x-y\|}{\diam K}
\smed\limits_{Q\in\, T(K)}\,\,\smed\limits_{|\xi|\le m}
(\diam K)^{|\xi|-m} \,|D^{\xi}P_{a_Q}(a_K)-D^{\xi}P_{a_K}(a_K)|
$$
holds. Here $C$ is a constant depending only on $n$ and $m$.
\end{lemma}
\par {\it Proof.} Notice that the function $F\in C^{\infty}(\RN\setminus E)$, see formula \rf{DEF-F}, so that, by the Lagrange theorem, for every $\alpha$, $|\alpha|=m$, there exists $z\in[x,y]$ such that
\bel{L-V}
|D^{\alpha}F(x)-D^{\alpha}F(y)|\le C\|x-y\|
\smed\limits_{|\beta|=m+1}\,\, |D^{\beta}F(z)|.
\ee
Since $x,y\in K^*$, the point $z\in K^*$ as well.
\par Since the polynomial $P_{a_K}\in\PMP$, for every multiindex $\beta$ of order $|\beta|=m+1$ we have $D^{\beta}P_{a_K}=0$ so that, by Lemma \reff{C-A},
\be
|D^{\beta}F(z)|&=&|D^{\beta}F(z)-D^{\beta}P_{a_K}(z)|\nn\\
&\le& C\,
\smed\limits_{Q\in T(K)}\,\,\smed\limits_{|\xi|\le m}
(\diam K)^{|\xi|-|\beta|} \,|D^{\xi}P_{a_Q}(a_K)-D^{\xi}P_{a_K}(a_K)|\nn\\&=&
C\,(\diam K)^{-1}
\smed\limits_{Q\in T(K)}\,\,\smed\limits_{|\xi|\le m}
(\diam K)^{|\xi|-m} \,|D^{\xi}P_{a_Q}(a_K)-D^{\xi}P_{a_K}(a_K)|.\nn
\ee
\par Combining this inequality with \rf{L-V} we obtain the statement of the lemma.\bx
\bigskip
\par We are in a position to prove inequality \rf{AIM} for arbitrary $y\in K$ and $x\in K^*$. By inequality \rf{W-P}, for every cube $Q\in T(K)$ and every $\xi, |\xi|\le m$,
$$
|D^{\xi}P_{a_Q}(a_K)-D^{\xi}P_{a_K}(a_K)|\le \lambda\,d(a_Q,a_K)(\diam K)^{m-|\xi|}
$$
so that, by Lemma \reff{L-V},
\be
I&:=&|D^{\alpha}F(x)-D^{\alpha}F(y)|\nn\\
&\le&
C\frac{\|x-y\|}{\diam K}
\smed\limits_{Q\in\, T(K)}\,\,\smed\limits_{|\xi|\le m}
(\diam K)^{|\xi|-m} \,|D^{\xi}P_{a_Q}(a_K)-D^{\xi}P_{a_K}(a_K)|\nn\\
&\le&
C\frac{\|x-y\|}{\diam K}
\smed\limits_{Q\in\, T(K)}\,\,\smed\limits_{|\xi|\le m}
(\diam K)^{|\xi|-m} \,(\lambda\,d(a_Q,a_K)(\diam K)^{m-|\xi|})\nn\\
&\le&
C\lambda\frac{\|x-y\|}{\diam K}
\smed\limits_{Q\in\, T(K)}\,d(a_Q,a_K).\nn
\ee
\par Notice that, by \rf{A-4} and \rf{DQ-E},
\be
\|a_Q-y\|&\le& \|a_Q-a_K\|+\|a_K-y\|\le 23\diam K+\diam K+\dist(K,E)\nn\\&\le& 23\cdot 4\dist(K,E)+4\dist(K,E)+\dist(K,E)\nn\\&=&97\dist(K,E)\le 97\,\|y-a_K\|.\nn
\ee
Hence, by part (b) of Claim \reff{VMT}, see \rf{O-1}, $d(a_Q,y)\le C\,d(y,a_K)$ so that
$$
d(a_Q,a_K)\le d(a_Q,y)+d(y,a_K)\le C\,d(y,a_K).
$$
\par This implies the following inequality
$$
I\le C\lambda\,(\#T(K))\,\frac{\|x-y\|}{\diam K}\,d(y,a_K).
$$
But, by part (2) of Lemma \reff{Wadd}, $\#T(K)\le N(n)$ so that
$$
I\le C\lambda\,\frac{\|x-y\|}{\diam K}\,d(y,a_K).
$$
\par Since $x,y\in K^*$, the distance
$\|x-y\|\le\diam K^*=\tfrac98\diam K.$ But
$$
\diam K\le 4\dist(K,E)\le 4\|y-a_K\|
$$
so that $\|x-y\|\le 5\|y-a_K\|$. Therefore, by part (b) of Claim \reff{VMT}, see \rf{O-2},
$$
\frac{\|x-y\|}{\|y-a_K\|}\,d(y,a_K)\le C\,d(x,y).
$$
\par On the other hand,
\be
\|y-a_K\|&\le& \diam K+\dist(a_K,K)=\diam K+\dist(K,E)\nn\\
&\le&
\diam K+4\diam K=5\diam K.\nn
\ee
Finally,
$$
I\le C\lambda\,\frac{\|x-y\|}{\diam K}\,d(y,a_K)\le
5\,C\lambda\,\frac{\|x-y\|}{\|y-a_K\|}\,d(y,a_K)
\le C\lambda\,d(x,y).
$$
\par The proof of Theorem \reff{CMD-E} is complete. \bx
\bigskip
\par {\bf 5.2 Extensions of $\LMP$-jets: a proof of Theorem \reff{JET-S}.}\medskip
\addtocontents{toc}{~~~~5.2 Extensions of $\LMP$-jets: a proof of Theorem \reff{JET-S}. \hfill \thepage\par}
\par Let $m\ge 1$ and let $J:E\to\PMRN$ be a mapping on $E$ which every point $x\in E$ assigns a polynomial $P_x=J(x)\in\PMRN$. We consider $J$ as a polynomial field
$$
J=\{P_x\in\PMRN: x\in E\}
$$
on $E$ and refer to $J$ as a {\it polynomial jet on $E$ of order $m-1$ or $(m-1)$-jet.}
\par Let
$$
JL^m_p(\RN)|_E=\{J:E\to\PMRN: \exists\,
F\in\LMP~~\text{such that}~~T^{m-1}_x[F]=P_x \,\,\forall x\in E\}
$$
be the space of traces to $E$ of all $(m-1)$-jets generated by $\LMP$-functions. We norm this space by the standard trace norm
\bel{J-NRM}
\|J\|_{JL^m_p(\RN)|_E}:=\inf \{\|F\|_{\LMP}:
F\in\LMP, T^{m-1}_x[F]=P_x~\text{for all}~ x\in E\}.
\ee
\par In these settings the result of Theorem \reff{JET-S} can  be reformulated as follows: for every $(m-1)$-jet
$$
J=\{P_x\in\PMRN: x\in E\}\in JL^m_p(\RN)|_E
$$
the following equivalence
$$
\|J\|_{JL^m_p(\RN)|_E}\sim \Nc_{m,p}(J)
$$
holds with constants depending only on $n,m,$ and $p$. Recall that the quantity $\Nc_{m,p}(J)$ is defined in the statement of Theorem \reff{JET-S}.
\medskip
\par Prove that
\bel{S-1}
\Nc_{m,p}(J)\le C(n,p,m)\,\|J\|_{JL^m_p(\RN)|_E}.
\ee
\par Let $F\in\LMP\cap \CMON$ be an arbitrary function such that $P_x=T^{m-1}_x[F]$ for all $x\in E$. Let $x\in\RN, y,z\in E$, and let
$$
K:=Q(x,2\|x-y\|+2\|x-z\|).
$$
By the \SP inequality \rf{SP-M}, for every multiindex $\beta, |\beta|\le m-1,$ we have
\be
|D^{\beta}P_y(y)-D^{\beta}P_z(y)|&\le& C\,\|y-z\|^{m-|\beta|}\left(\frac{1}{|Q_{yz}|}
\intl_{Q_{yz}}(\nabla ^mF(u))^qdu\right)^{\frac{1}{q}}\nn\\&\le&
C\,\|y-z\|^{m-|\beta|-\tfrac{n}{q}}
\left(\,\intl_{Q_{yz}}
(\nabla ^mF(u))^qdu\right)^{\frac{1}{q}}.\nn
\ee
Since $Q_{yz}\subset K$ and $m-|\beta|-\tfrac{n}{q}>0$ provided $|\beta|\le m-1$, we obtain
\be
|D^{\beta}P_y(y)-D^{\beta}P_z(y)|&\le&
C(\|x-y\|+\|x-z\|)^{m-|\beta|-\tfrac{n}{q}}
\left(\,\intl_{K}
(\nabla ^mF(u))^qdu\right)^{\frac{1}{q}}\nn\\
&\le&
C(\,\|x-y\|^{m-|\beta|}+\|x-z\|^{m-|\beta|})\,
|K|^{-\tfrac{n}{q}} \left(\,\intl_{K}
(\nabla^m F(u))^qdu\right)^{\frac{1}{q}}.\nn
\ee
Hence
\bel{DJ-A}
\frac{|D^{\beta}P_y(y)-D^{\beta}P_z(y)|}
{\|x-y\|^{m-|\beta|}+\|x-z\|^{m-|\beta|}}
\le
C\,\left(\frac{1}{|K|}\,\intl_{K}
(\nabla ^mF(u))^qdu\right)^{\frac{1}{q}}.
\ee
\par Let us introduce a sharp maximal function  $J^{\#}_E:\RN\to\R_+$ for the jet $J=\{P_x: x\in E\}$ by letting
\bel{J-SH}
J^{\#}_E(x):=
\sbig_{|\beta|\le m-1}\,\,
\sup_{y,z\in E}\,\, \frac{|D^{\beta}P_y(y)-D^{\beta}P_z(y)|}{\|x-y\|^{m-|\beta|}+
\|x-z\|^{m-|\beta|}}\,\,,~~~~~x\in\RN.
\ee
Clearly, by \rf{IDF},
\bel{I-J}
\|J^{\#}_E\|_{\LPRN}\sim \Nc_{m,p}(J)\,.
\ee
\smallskip
\par Prove that $\|J^{\#}_E\|_{\LPRN}\le C\|F\|_{\LMP}$. By inequality \rf{DJ-A},
$$
\sup_{y,z\in E}\,\, \frac{|D^{\beta}P_y(y)-D^{\beta}P_z(y)|}{\|x-y\|^{m-|\beta|}+
\|x-z\|^{m-|\beta|}}\le C\, \left\{\Mc[(\nabla^m F(u))^q](x)\right\}^{\frac1q}
$$
so that
$$
J^{\#}_E(x)\le C\, \left\{\Mc[(\nabla^m F)^q](x) \right\}^{\frac1q},~~~~x\in\RN.
$$
Since $q<p$, by the Hardy-Littlewood maximal theorem,
$$
\|J^{\#}_E\|_{\LPRN}\le C\, \|\left
\{\Mc[(\nabla^m F))^q]\right\}^{\frac1q}\|_{\LPRN}
\le
C\, \|\nabla^m F\|_{\LPRN}= C\,\|F\|_{\LMP}
$$
so that, by \rf{I-J}, $\Nc_{m,p}(J)\le C\,\|F\|_{\LMP}$. Taking the infimum in this inequality over all functions $F\in\LMP\cap \CMON$ such that
$T^{m-1}_x[F]=P_x$ for every $x\in E$, we obtain the required inequality \rf{S-1}.
\medskip
\par Prove that
\bel{S-2}
\|J\|_{JL^m_p(\RN)|_E}\le C(n,p,m)\,\Nc_{m,p}(J).
\ee
\par Let $q:=(n+p)/2$ and $\theta:=(q+p)/2$. Let $h_1:=\Mc[(J^{\#}_E)^\theta]^{\frac{1}{\theta}}$. Clearly, $J^{\#}_E\le h_1$ a.e. on $\RN$.
\par Let $x,y\in E$ and let $u\in Q_{xy}=Q(x,\|x-y\|)$. Then
$$
\|u-x\|,\,\|u-y\|\le 2\|x-y\|
$$
so that, by \rf{J-SH},
$$
|D^{\beta}P_x(x)-D^{\beta}P_y(x)|\le (\|u-x\|^{m-|\beta|}+\|u-y\|^{m-|\beta|})\, J^{\#}_E(u)
\le C\,\|x-y\|^{m-|\beta|}J^{\#}_E(u)
$$
for every $\beta, |\beta|\le m-1,$ and every $u\in Q_{xy}$. Integrating this inequality over the cube $Q_{xy}$ with respect to $u$, we obtain
\be
|D^{\beta}P_x(x)-D^{\beta}P_y(x)|
&\le& C\,\|x-y\|^{m-|\beta|}\left(\frac{1}{|Q_{xy}|}\,
\intl_{Q_{xy}}
J^{\#}_E(u)\,du\right)\nn\\&\le& C\,\|x-y\|^{m-|\beta|}\left(\frac{1}{|Q_{xy}|}\,
\intl_{Q_{xy}}
h_1(u)\,du\right)\nn
\ee
so that
$$
|D^{\beta}P_x(x)-D^{\beta}P_y(x)|
\le C\,\|x-y\|^{m-|\beta|}\left(\frac{1}{|Q_{xy}|}\,
\intl_{Q_{xy}}
h^q_1(u)\,du\right)^{\frac1q}.
$$
Hence
$$
\frac{|D^{\beta}P_x(x)-D^{\beta}P_y(x)|}{\|x-y\|^{m-|\beta|-1}}
\le C\,\|x-y\|\left(\frac{1}{|Q_{xy}|}\,
\intl_{Q_{xy}}
h^q_1(u)\,du\right)^{\frac1q}
$$
proving that
\bel{D-BT}
\frac{|D^{\beta}P_x(x)-D^{\beta}P_y(x)|}
{\|x-y\|^{m-|\beta|-1}}
\le \delta_q(x,y:h_2),~~~~x,y\in E,
\ee
where $h_2:=C\,h_1$. See \rf{DQ}.
\par Similar to \rf{H-INA}, \rf{H-A1} and \rf{H-LP} we show that
\bel{P-H}
h^q_2\in A_1(\RN),~~\|h^q_2\|_{A_1}\le\eta(n,p,q),
\ee
and
\bel{H-JP}
\|h_2\|_{\LPRN}\le C\,\|J^{\#}_E\|_{\LPRN}.
\ee
\par Let $d=d_q(x,y:h_2)$, $x,y\in\RN,$ be the geodesic metric associated with $\delta_q(h_2)$, see \rf{GD-M}. By \rf{P-H} and \rf{H-JP}, the metric $d\in\DPQ$, see \rf{D-LIP}. Hence, by Theorem \reff{MAIN-MT},
$$
\dq(x,y:h_2)\le C\, d_q(x,y:h_2),~~~~x,y\in\RN,
$$
where $C=C(n,p)$.
\par This inequality and \rf{D-BT} imply that for every $x,y\in E$ and every $\beta, |\beta|\le m-1$,
$$
\frac{|D^{\beta}P_x(x)-D^{\beta}P_y(x)|}
{\|x-y\|^{m-|\beta|-1}}
\le C\,d_q(x,y:h_2)=C\,d(x,y).
$$
Hence
$$
S[J]:=\sbig_{|\beta|\le m-1}\,\,
\sup_{x,y\in E,\,x\ne y}\,\, \frac{|D^{\beta}P_y(y)-D^{\beta}P_z(y)|}
{\|x-y\|^{m-|\beta|-1}\,d(x,y)}\le C(n,m,p), ~~~~x,y\in\RN,
$$
so that, by Theorem \reff{CMD-E}, the jet $J\in JC^{m-1,(d)}(\RN)|_E$.
\par Thus there exists a function $F\in C^{m-1,(d)}(\RN)$ such that $T^{m-1}_x[F]=P_x$ for every $x\in E$. Furthermore
$$
\|F\|_{C^{m-1,(d)}(\RN)}\le C\, S[J]\le C(n,m,p),
$$
so that for every $\alpha, |\alpha|=m-1,$ and every $x,y\in\RN$
$$
|D^\alpha F(x)-D^\alpha F(y)|\le C\,d_q(x,y:h_2)=d_q(x,y:h_3)
$$
where $h_3=C\,h_2$. Since $h^q_2\in A_1(\RN)$ and $\|h^q_2\|_{A_1}\le\eta(n,p,q)$, the same is true for the function $h_3$ as well.
\par Hence, by Theorem \reff{SM-MAIN}, the function $F\in\LMP$ and
$$
\|F\|_{\LMP}\le C\,\|h_3\|_{\LMP}\le C\,\|h_2\|_{\LMP}.
$$
so that, by \rf{H-JP} and \rf{I-J},
$$
\|F\|_{\LMP}\le C\,\|J^{\#}_E\|_{\LPRN}\le C\, \Nc_{m,p}(J).
$$
Finally,
$$
\|J\|_{JL^m_p(\RN)|_E}\le \|F\|_{\LMP}\le C(n,m,p)\,\Nc_{m,p}(J)
$$
proving inequality \rf{S-2}.
\par Theorem \reff{JET-S} is completely proved.\bx
\bigskip
\par {\bf 5.3 A variational criterion for the traces of  $\LMP$-jets.}\medskip
\addtocontents{toc}{~~~~5.3 A variational criterion for the traces of $\LMP$-jets. \hfill \thepage\par}
\par In \cite{Sh2} we have proved a criterion which provides a characterization of the trace space $\LOP|_E$
it terms of certain local oscillations of functions on subsets of the set $E$. The next theorem generalizes this result to the case of jet-spaces generated by $L^m_p(\RN)-$functions.
\begin{theorem}\lbl{EX-TK} Let $m\in\N$ and let $p\in (n,\infty)$. Let $J=\{P_x\in\Pc_{m-1}:x\in E\}$ be a
polynomial field on a closed set $E\subset\RN$.
\par There exists a $C^{m-1}$-function $F\in\LMP$ such that $T_{x}^{m-1}(F)=P_{x}$ for every $x\in E$ if and only if there exists a constant $\lambda>0$ such
that for every finite family $\{Q_i:i=1,...,k\}$ of
disjoint cubes in $\RN$, every $x_i,y_i\in (\gamma Q_i)\cap E,$ and every multiindex $\beta, |\beta|\le m-1,$ the
following inequality
\bel{N-P}
\sbig_{i=1}^k\,\,\frac{|D^\beta P_{x_i}(x_i)
-D^\beta P_{y_i}(x_i)|^p}
{(\diam Q_i)^{(m-|\beta|)p-n}} \le\lambda
\ee
holds. Here $\gamma>1$ is an absolute constant.
\par Furthermore,
$$
\|J\|_{JL^m_p(\RN)|_E}\sim \inf\lambda^{\frac1p}.
$$
See \rf{J-NRM}. The constants in this equivalence depend only on $n,m$ and $p$.
\end{theorem}
\bigskip
\par {\it Proof.} {\it (Necessity.)} Let
$$
J=\{P_x\in\PMRN:x\in E\}
$$
be a polynomial field defined on $E$. Let $\gamma>1$ be a constant and let $F\in\LMP\cap \CMON$ be an arbitrary function such that $P_x=T^{m-1}_x[F]$ for all $x\in E$. Let $Q$ be a cube in $\RN$ and let $x,y\in K:=\gamma Q$. Then, by the \SP inequality \rf{SP-M}, for
every multiindex $\beta$, $|\beta|\le m-1,$ we have
\be
|D^{\beta}P_x(x)-D^{\beta}P_y(x)|&\le& C\,\|x-y\|^{m-|\beta|}\left(\frac{1}{|Q_{xy}|}
\intl_{Q_{xy}}(\nabla ^mF(u))^qdu\right)^{\frac{1}{q}}\nn\\
&\le&
C\,\|x-y\|^{m-|\beta|-\tfrac{n}{p}}
\left(\,\intl_{Q_{xy}}(\nabla ^mF(u))^qdu\right)^{\frac{1}{q}}\nn.
\ee
Clearly, $Q_{xy}\subset 2K$ and $m-|\beta|-\tfrac{n}{p}>0$ whenever $n<p$ and $|\beta|\le m-1$ so that
$$
|D^{\beta}P_x(x)-D^{\beta}P_y(x)|\le C\,
(\diam Q)^{m-|\beta|} \left(\frac{1}{|2K|}
\intl_{2K}(\nabla ^mF(u))^qdu\right)^{\frac{1}{q}}.
$$
Therefore for every $z\in Q$ we have
$$
|D^{\beta}P_x(x)-D^{\beta}P_y(x)|^p\le C\,
(\diam Q)^{p(m-|\beta|)} \left(\Mc[\nabla ^m F^q](z)\right)^{\frac{p}{q}}.
$$
Integrating this inequality over cube $Q$
(with respect to $z$) we obtain
$$
\frac{|D^{\beta}P_x(x)-D^{\beta}P_y(x)|^p}
{(\diam Q)^{(m-|\beta|)p-n}}\le C\,\intl_{Q}
\left(\Mc[\nabla ^m F^q](z)\right)^{\frac{p}{q}}\,dz.
$$
\par Hence,
\be
I_\beta:=\sbig_{i=1}^m\,\,\frac{|D^\beta P_{x_i}(x_i)
-D^\beta P_{y_i}(x_i)|^p}
{(\diam Q_i)^{(m-|\beta|)p-n}}& \le&
C\,\sbig_{i=1}^m\,
\intl_{Q_i}
\left(\Mc[\nabla ^m F^q](z)\right)^{\frac{p}{q}}\,dz\nn\\
&\le& C\,\intl_{\RN}
\left(\Mc[\nabla ^m F^q](z)\right)^{\frac{p}{q}}\,dz\nn
\ee
so that, by the Hardy-Littlewood maximal theorem,
\bel{IB-F}
I_\beta\le C\,\intl_{\RN}
(\nabla ^m F)^p(z)\,dz\,\le C\,\|F\|^p_{\LMP}\,.
\ee
This proves \rf{N-P} with $\lambda=C\,\|F\|^p_{\LMP}.$ Furthermore, taking the infimum in \rf{IB-F} over all functions $F\in\LMP$ such that $T^{m-1}_x[F]=P_x$ on $E$ we obtain that
$$
I_\beta\le C\,\|J\|_{JL^m_p(\RN)|_E}
$$
proving the necessity.
\bigskip
\par {\it (Sufficiency).} Let $\gamma:=\,10^4$. Let  $J=\{P_x\in\PMRN:x\in E\}$ be a polynomial field on $E$. Suppose that there exists a constant $\lambda>0$ such that for every family $\{Q_i:i=1,...,k\}$ of pairwise disjoint cubes in $\RN$, every $x_i,y_i\in (\gamma Q_i)\cap E,$ and every multiindex $\beta, |\beta|\le m-1,$ inequality \rf{N-P} is satisfied. Here $\gamma$ is a certain absolute constant which we determine below.
\smallskip
\par Let us prove that under these conditions the Whitney extension $F:\RN\to\R$ of the jet $J$ defined by the formula \rf{DEF-F} has the following properties:\smallskip
\par (i). $F\in\CMON$ and $T_x^{m-1}[F]=P_x$ for every $x\in E$\,;\smallskip
\par (ii). $F\in\LMP$ and
\bel{P-2}
\|F\|_{\LMP}\le C\,\lambda^{\frac1p}.
\ee
\par Prove (i). For every multiindex $\beta, |\beta|\le m-1,$ and every $x,y\in E$, by \rf{N-P},
$$
\frac{|D^\beta P_{x}(x)-D^\beta P_{y}(x)|^p}
{(\diam Q_{xy})^{(m-|\beta|)p-n}} \le\lambda.
$$
Recall that $Q_{xy}=Q(x,\|x-y\|)$ so that $\diam Q_{xy}=2\|x-y\|$. Hence
$$
|D^\beta P_{x}(x)-D^\beta P_{y}(x)|
\le C\,\lambda^{\frac1p}\|x-y\|^{m-|\beta|-1}\cdot
\|x-y\|^{1-\frac{n}{p}}.
$$
But $n<p$ so that
$$
D^\beta P_{x}(x)-D^\beta P_{y}(x)
=o(\|x-y\|^{m-|\beta|-1})~~\text{as}~~y\to x,~y\in E.
$$
Thus the jet $J=\{P_x\in\PMRN:x\in E\}$ satisfies the hypothesis of the Whitney extension theorem \cite{W1}. This  theorem implies the statement (i).
\smallskip
\par Prove (ii). Let us fix a multiindex  $\beta$ of order $|\beta|=m-1$ and prove that $D^\beta F\in \LOP$ and
$$
\|D^\beta F\|_{\LOP}\le C\,\lambda^{\frac1p}.
$$
\par We will make use of a characterization of Sobolev spaces which follows from results proven in \cite{Br}: Let $p>n$ and let $G\in C(\RN)$. Suppose that there exists a constant $\tau>0$ such that the following inequality
\bel{G-C}
\sbig_{i=1}^k\,\frac{|G(x_i)-G(x_i)|^p}{(\diam Q_i)^{p-n}}
\le \tau
\ee
holds for every finite family $\{Q_i: i=1,...,k\}$ of pairwise disjoint equal cubes and all $x_i,y_i\in Q_i$. Then $G\in\LOP$ and $\|G\|_{\LOP}\le C(n,p)\,\tau^{\frac1p}$.
\medskip
\par We are also needed the following combinatorial
\begin{theorem}(\cite{BrK,Dol})\lbl{TFM} Let $\Ac=\{Q\}$ be a collection of cubes in $\RN$ with covering
multiplicity $ \MP(\Ac)<\infty$. Then $\Ac$ can be
partitioned into at most $N=2^{n-1}(\MP(\Ac)-1)+1$ families
of disjoint cubes.
\end{theorem}
\par Recall that {\it covering multiplicity} $\MP(\Ac)$ of a family of cubes $\Ac$ is the minimal positive integer $M$ such that every point $x\in\RN$ is covered by at most $M$ cubes from $\Ac$.
\begin{lemma}\lbl{L-BG} Let $\Qc=\{Q_1,...,Q_k\}$ be a family of pairwise disjoint equal cubes in $\RN$ such that
$$
\dist(c_{Q_i},E)\le 40\diam Q_i,~~~i=1,...,k.
$$
Then for every $x_i,y_i\in Q_i$ the following inequality
$$
\sbig_{i=1}^k\,
\frac{|D^\beta F(x_i)-D^\beta F(y_i)|^p}
{(\diam Q_i)^{p-n}}
\le C\,\lambda
$$
holds.
\end{lemma}
\par {\it Proof.} Fix $i\in\{1,...,k\}$. Let $Q=Q_i\in\Qc$ so that
\bel{D80}
\dist(c_{Q},E)\le 40\diam Q\,.
\ee
\par Let $p_Q\in E$ be a point nearest to $Q$ on $E$. Then
$$
\|c_Q-p_Q\|\le\diam Q+40\diam Q=41\diam Q=82 r_Q.
$$
Hence $p_Q\in 82 Q$.
\par Let $x_Q=x_i$ and $y_Q=y_i$. (Recall that $Q=Q_i$ for some $i\in\{1,...,k\}$.) Then
\be
I_Q&:=&\frac{|D^\beta F(x_Q)-D^\beta F(y_Q)|^p}
{(\diam Q)^{p-n}}\nn\\&\le& 2^p\left\{
\frac{|D^\beta F(x_Q)-D^\beta F(p_Q)|^p}{(\diam Q)^{p-n}}+\frac{|D^\beta F(p_Q)-D^\beta F(y_Q)|^p}
{(\diam Q)^{p-n}}\right\}\nn
\ee
so that
\be
I&:=&\smed_{Q\in\Qc}\,\,I_Q\le C\,\left\{
\sbig_{Q\in\Qc}\,
\frac{|D^\beta F(x_Q)-D^\beta F(p_Q)|^p}{(\diam Q)^{p-n}}\right.\nn\\
&+&\left.\sbig_{Q\in\Qc}\,
\frac{|D^\beta F(p_Q)-D^\beta F(y_Q)|^p}
{(\diam Q)^{p-n}}\right\}=C\,\{I_1+I_2\}\,.\nn
\ee
\par Recall that $p_Q\in E$ so that
$D^\beta F(p_Q)=D^\beta P_{p_Q}(p_Q)$.
\par Now suppose that $x_Q\in E$. Then $D^\beta F(x_Q)=D^\beta P_{x_Q}(x_Q)$. Furthermore, since
$P_{x_Q}\in \PMRN$ and $|\beta|=m-1$, the function
$D^\beta P_{x_Q}$ is constant so that $D^\beta F(x_Q)=D^\beta P_{x_Q}(p_Q)$. Hence
$$
\frac{|D^\beta F(x_Q)-D^\beta F(p_Q)|^p}
{(\diam Q)^{p-n}}
=\frac{|D^\beta P_{x_Q}(p_Q)-D^\beta P_{p_Q}(p_Q)|^p}
{(\diam Q)^{p-n}}
$$
so that, by assumption \rf{N-P} (with $|\beta|=m-1$), we have
$$
I_1:=\sbig_{Q\in\Qc}\,
\frac{|D^\beta F(x_Q)-D^\beta F(p_Q)|^p}{(\diam Q)^{p-n}}\le \lambda
$$
provided $x_Q\in E$ for every $Q\in\Qc$.
\par Thus later on we can assume that $x_Q\in\RN\setminus E$ for all $Q\in\Qc$.
\par Let $K_Q\in W_E$ be {\it a Whitney cube which contains $x_Q$.} Recall that given $H\in W_E$ by $a_H$ we denote a point nearest to $H$ on $E$. Also by $T(K_Q)$ we denote the family of Whitney's cubes intersecting $K_Q$. See \rf{TK}.
\par Let
$$
S(Q):=|D^{\beta}F(x_Q)-D^{\beta}F(p_Q)|
$$
and let
$$
V(Q):=|D^{\beta}P_{p_Q}(p_Q)-D^{\beta}P_{a_{K_Q}}(p_Q)|.
$$
Given $H\in T(K_Q)$ and a multiindex $\xi$ with $|\xi|\le m-1$ let
$$
L(\xi:H,Q):=|D^{\xi}P_{a_{H}}(a_{K_Q})
-D^{\xi}P_{a_{K_Q}}(a_{K_Q})|.
$$
\par Then, by Lemma \reff{C-2},
$$
S(Q)\le
C\,\left\{V(Q)+
\sbig\limits_{H\in T(K_Q)}\,\,\sbig\limits_{|\xi|\le m-1}\frac{L(\xi:H,Q)}{(\diam K_Q)^{m-1-|\xi|}}\right\}.
$$
Since $\#T(K_Q)\le N(n),$ see Lemma \reff{Wadd}, we have
$$
\frac{S(Q)^p}
{(\diam Q)^{p-n}}\le
\frac{C\,V(Q)^p}
{(\diam Q)^{p-n}}+
C\sbig\limits_{H\in T(K_Q)}\,\,
\sbig\limits_{|\xi|\le m-1}
\,\left(\frac{\diam K_Q}{\diam Q}\right)^{p-n}\frac{L(\xi:H,Q)^p}
{(\diam K_Q)^{(m-|\xi|)p-n}}.
$$
\par Prove that $a_{K_Q}\in \gamma Q$. In fact, since $x_Q\in K_Q\cap Q$, we have
\be
\diam K_Q&\le& 4\dist(K_Q,E)\le 4\dist(x_Q,E)\nn\\
&\le& 4(\|x_Q-c_Q\|+\dist(c_Q,E))\le 4\diam Q+4\dist(c_Q,E)\nn
\ee
so that, by \rf{D80},
$$
\diam K_Q\le 4\diam Q+4\cdot 40\diam Q=164\diam Q.
$$
In particular, by this inequality,
$$
\frac{S(Q)^p}
{(\diam Q)^{p-n}}\le
C\,\left\{\frac{V(Q)^p}
{(\diam Q)^{p-n}}+
\sbig\limits_{H\in T(K_Q)}\,\,
\sbig\limits_{|\xi|\le m-1}
\,\frac{L(\xi:H,Q)^p}
{(\diam K_Q)^{(m-|\xi|)p-n}}\right\}.
$$
\par Since $K_Q\cap Q\ne\emp$, we have
\be
\|c_Q-a_{K_Q}\|&\le& \dist(a_{K_Q},K_Q)+\diam K_Q+\diam Q\nn\\
&=&
\dist(K_Q,E)+\diam K_Q+\diam Q\nn\\
&\le& 4\diam K_Q+\diam K_Q+\diam Q\nn\\
&\le&
5\cdot 164\diam Q+\diam Q\le (\gamma/2)\diam Q.\nn
\ee
(Recall that $\gamma:=10^4.$) Hence $a_{K_Q}\in \gamma Q$. 
\par Prove that $a_H\in\gamma K_Q$ whenever $H\in T(K_Q)$. Since $H\cap K_Q\ne\emp$ and $\diam H\le 4\diam K_Q$, see Lemma \reff{Wadd}, we have
\be
\|c_{K_Q}-a_{H}\|&\le& \dist(a_{H},E)+\diam H+\diam K_Q\nn\\
&=&
\dist(H,E)+\diam H+\diam K_Q\nn\\
&\le& 4\diam H+\diam H+\diam K_Q\nn\\
&\le&
5\cdot 4\diam K_Q+\diam K_Q=21\diam K_Q.\nn
\ee
Hence
\bel{A-TK}
a_H\in 42 K_Q~~~\text{for every}~~~H\in T(K_Q)
\ee
proving that $a_H\subset\gamma K_Q$.
\par By $H_Q$ we denote a cube $H\in T(K_Q)$ for which the quantity
$$
\sbig\limits_{|\xi|\le m-1}
\,\frac{L(\xi:H,Q)^p}
{(\diam K_Q)^{(m-|\xi|)p-n}}
$$
is maximal on $T(K_Q)$. Since $\#T(K_Q)\le N(n),$ we obtain the following inequality
$$
\frac{S(Q)^p}
{(\diam Q)^{p-n}}\le
C\,\left\{\frac{V(Q)^p}
{(\diam Q)^{p-n}}+
\sbig\limits_{|\xi|\le m-1}
\,\frac{L(\xi:H_Q,Q)^p}
{(\diam K_Q)^{(m-|\xi|)p-n}}\right\}.
$$
Hence
\be
I_2&:=&\sbig_{Q\in\Qc}\frac{S(Q)^p}
{(\diam Q)^{p-n}}\nn\\
&\le&
C\,\left\{\sbig_{Q\in\Qc}\frac{V(Q)^p}
{(\diam Q)^{p-n}}+
\sbig\limits_{|\xi|\le m-1}\sbig_{Q\in\Qc}
\frac{L(\xi:H_Q,Q)^p}
{(\diam K_Q)^{(m-|\xi|)p-n}}\right\}\nn\\
&=&
C\,\{I_2^{(1)}+I_2^{(2)}\}.
\nn
\ee
\par Since $p_Q\in 82 Q\subset\gamma Q$ and $a_{K_Q}\in\gamma Q$ we can apply assumption \rf{N-P} to the family of cubes $\Qc$. By this assumption $I_2^{(1)}\le \lambda$.
\par Let $\Kc:=\{K_Q: Q\in\Qc\}$. We know that the cubes of this family are non-overlapping, but, in general, they are not disjoint so that we can not apply the assumption of the theorem to $\Kc$. Nevertheless, the family $\Kc$ as a subfamily of $W_E$ has covering multiplicity bounded by a constant $N(n)$. By Theorem \reff{TFM}, {\it it can be partitioned into at most $M(n)$ families of pairwise disjoint cubes} so that, without loss of generality, we may assume {\it $\Kc$ itself consists of pairwise disjoint cubes.}
\par Since $a_{H_Q}, a_{K_Q}\in\gamma K_Q$ for every $Q\in\Qc$, by assumption \rf{N-P},
$$
I_2^{(2)}:=
\sbig\limits_{|\xi|\le m-1}\,\,
\sbig_{Q\in\Qc}\,
\frac{|D^{\xi}P_{a_{H_Q}}(a_{K_Q})
-D^{\xi}P_{a_{K_Q}}(a_{K_Q})|^p}
{(\diam K_Q)^{(m-|\xi|)p-n}}\le
\smed\limits_{|\xi|\le m-1}\lambda\le C\lambda
$$
proving the lemma.\bx
\smallskip
\begin{lemma}\lbl{L-SMC} Let $\Qc=\{Q_1,...,Q_k\}$ be a family of pairwise disjoint equal cubes in $\RN$ such that
\bel{D-S80}
\diam Q_i<\tfrac{1}{40}\dist(c_{Q_i},E),~~~i=1,...,k.
\ee
Then
$$
\sbig_{i=1}^k\,
\frac{|D^\beta F(x_i)-D^\beta F(y_i)|^p}
{(\diam Q_i)^{p-n}}
\le C\,\lambda
$$
for every choice of points $x_i,y_i\in Q_i$.
\end{lemma}
\par {\it Proof.} For every cube $Q\in\Qc$, by \rf{D-S80}, $\dist(c_{Q},E)> 40\diam Q>0$ so that $Q\subset \RN\setminus E$. Let $K_Q$ be a Whitney cube which contains $c_Q$. Prove that $Q\subset K^*_Q=\tfrac98 K_Q.$
\par In fact,
\be
\diam Q&<&\tfrac1{40}\dist(c_{Q},E)\le \tfrac1{40}\{\diam K_Q+\dist(K_{Q},E)\}\nn\\&\le& \tfrac1{40}
\{\diam K_Q+4\diam K_Q\}=\tfrac1{8}\diam K_Q.\nn
\ee
Hence for every $z\in Q$ we have
\be
\|z-c_{K_Q}\|&\le& \|z-c_Q\|+\|c_Q-c_{K_Q}\|\le \tfrac12\diam Q+\tfrac12\diam K_Q\nn\\
&\le& \tfrac12\cdot\tfrac1{8}\diam K_Q+\tfrac12\diam K_Q=\left(\tfrac1{8}+1\right) r_{K_Q}\nn
\ee
so that $Q\subset\tfrac98 K_Q=K^*_Q$. Hence the points
$$
x_Q=x_i,~y_Q=y_i\in K^*_Q
$$
provided $Q=Q_i$.
\par Let $K\in W_E$ and let
$$
\Qc(K):=\{Q\in\Qc:K=K_Q\}=\{Q\in\Qc: c_Q\in K\}.
$$
By $\Kc$ we denote a family of Whitney's cubes $K$ for which $\Qc(K)\ne\emp$.
\par Then for each $K\in \Kc$, by Lemma \reff{C-4},
$$
|D^{\beta}F(x_Q)-D^{\beta}F(y_Q)|\le C
\frac{\|x_Q-y_Q\|}{\diam K}
\sbig\limits_{H\in\, T(K)}\,\,\sbig\limits_{|\xi|\le m}
 \,\frac{|D^{\xi}P_{a_H}(a_K)-D^{\xi}P_{a_K}(a_K)|}{(\diam K)^{m-|\xi|}}.
$$
\par By \rf{A-TK},
\bel{IN-A}
a_K,a_H\in 42 K~~~\text{for every}~~~H\in T(K).
\ee
\par Now we have
\be
I_K&:=&
\sbig_{Q\in\Qc(K)}\,
\frac{|D^\beta F(x_Q)-D^\beta F(y_Q)|^p}
{(\diam Q)^{p-n}}\nn\\
&\le&
C\left\{\sbig_{Q\in\Qc(K)}
\left(\frac{\|x_Q-y_Q\|}{\diam Q}\right)^p|Q|\right\}
\left\{\sbig\limits_{H\in T(K)}\,\,
\sbig\limits_{|\xi|\le m}
\frac{|D^{\xi}P_{a_H}(a_K)-D^{\xi}P_{a_K}(a_K)|}{(\diam K)^{m-|\xi|}}\right\}^p\nn\\
&\le&
C\,\left\{\sbig_{Q\in\Qc(K)}\,
|Q|\right\}
\left\{\sbig\limits_{H\in\, T(K)}\,\,
\sbig\limits_{|\xi|\le m}
\,\frac{|D^{\xi}P_{a_H}(a_K)-D^{\xi}P_{a_K}(a_K)|}{(\diam K)^{m-|\xi|}}\right\}^p\nn\\
&\le&
C\,|K|\left\{\sbig\limits_{H\in\, T(K)}\,\,
\sbig\limits_{|\xi|\le m}
\,\frac{|D^{\xi}P_{a_H}(a_K)-D^{\xi}P_{a_K}(a_K)|}{(\diam K)^{m-|\xi|}}\right\}^p.\nn
\ee
\par Since $\#T(K)\le N(n)$, see Lemma \reff{Wadd}, we obtain
$$
I_K\le
C\,\sbig\limits_{H\in\, T(K)}\,\,
\sbig\limits_{|\xi|\le m}
\,\frac{|D^{\xi}P_{a_H}(a_K)-D^{\xi}P_{a_K}(a_K)|^p}
{(\diam K)^{(m-|\xi|)p-n}}.
$$
Let $\tH\in T(K)$ be a cube such that the quantity
$$
\sbig\limits_{|\xi|\le m}
\,\frac{|D^{\xi}P_{a_H}(a_K)-D^{\xi}P_{a_K}(a_K)|^p}
{(\diam K)^{(m-|\xi|)p-n}}
$$
takes the maximal value on $T(K)$. Then
$$
I_K\le
C\,\sbig\limits_{|\xi|\le m}
\,\frac{|D^{\xi}P_{a_{\tH}}(a_K)-D^{\xi}P_{a_K}(a_K)|^p}
{(\diam K)^{(m-|\xi|)p-n}}.
$$
Hence
\be
I&:=&\sbig_{i=1}^k\,
\frac{|D^\beta F(x_i)-D^\beta F(y_i)|^p}
{(\diam Q_i)^{p-n}}\le\smed_{K\in \Kc(\Qc)} I_K\nn\\
&\le&
C\,\sbig\limits_{|\xi|\le m}\,\,\sbig_{K\in \Kc(\Qc)}
\,\frac{|D^{\xi}P_{a_{\tH}}(a_K)-D^{\xi}P_{a_K}(a_K)|^p}
{(\diam K)^{(m-|\xi|)p-n}}.\nn
\ee
\par As in the proof of the previous lemma, using Theorem \reff{TFM} we can assume that the cubes of the family $\Kc(\Qc)$ are pairwise disjoint. This and inclusions
\rf{IN-A} enable us to apply assumption \rf{N-P} to the last sum of the above inequality. By this assumption
$$
I\le C\,\smed\limits_{|\xi|\le m}\,\lambda\le C\,\lambda
$$
proving the lemma.\bx
\par Combining Lemma \reff{L-BG} and Lemma \reff{L-SMC} with the criterion \rf{G-C} we conclude that $F$ is a $C^{m-1}$-smooth function such that for every multiindex $\beta$ of order $m-1$ the function $D^\beta F\in\LOP$ and $\|D^\beta F\|_{\LOP}\le C\,\lambda^{\frac1p}.$ Hence $F\in\LMP$ and
$\| F\|_{\LMP}\le C\,\lambda^{\frac1p}$.
\par The proof of Theorem \reff{EX-TK} is complete.\bx

\end{document}